\documentclass[12pt,a4paper,reqno]{amsart}
\usepackage{amsmath}
\usepackage{amsfonts}
\usepackage{amssymb}
\numberwithin{equation}{section}

\theoremstyle{plain}

\newtheorem{theorem}[subsection]{Theorem}
\newtheorem{proposition}[subsection]{Proposition}
\newtheorem{lemma}[subsection]{Lemma}
\newtheorem{corollary}[subsection]{Corollary}

\theoremstyle{definition}

\newtheorem{definition}[subsection]{Definition}

\renewcommand{\leq}{\leqslant}
\renewcommand{\geq}{\geqslant}

\newsavebox{\proofbox}
\savebox{\proofbox}{\begin{picture}(7,7)%
  \put(0,0){\framebox(7,7){}}\end{picture}}
\def\boxeq{\tag*{\usebox{\proofbox}}}

\newcommand{\md}[1]{\ensuremath{(\mbox{mod}\, #1)}}

\def\endproof{\hfill{\usebox{\proofbox}}}
\def\E{\mathbb{E}}
\def\N{\mathbb{N}}
\def\F{\mathbb{F}}
\def\Z{\mathbb{Z}}
\def\R{\mathbb{R}}
\def\T{\mathbb{T}}
\def\C{\mathbb{C}}
\def\P{\mathbb{P}}

\def\B{\mathcal{B}}

\def\rank{\operatorname{rank}}
\def\eps{\varepsilon}

\def\f{{\mathbf f}}
\def\F{{\mathbf F}}
\def\g{{\mathbf g}}
\def\lin{\operatorname{lin}}

\begin{document}

\title[Length 4 progressions in finite field geometries]{New bounds for Szemer\'edi's theorem, I: Progressions of length 4 in finite field geometries}

\author{Ben Green}
\address{Centre for Mathematical Sciences, Wilberforce Road, Cambridge CB3 0WA.}
\email{b.j.green@dpmms.cam.ac.uk}

\author{Terence Tao}
\address{Department of Mathematics, UCLA, Los Angeles CA 90095-1555, USA.
}
\email{tao@math.ucla.edu}

\thanks{The first author is a Clay Research Fellow, and is pleased to acknowledge the support of the Clay Mathematics Institute. Some of this work was carried out while he was on a long-term visit to MIT. The second author is supported by a grant from the Packard Foundation.}

\begin{abstract} 
 Let $k \geq 3$ be an integer, and let $G$ be a finite abelian group with $|G| = N$, where $(N, (k-1)!) = 1$. We write $r_k(G)$ for the largest cardinality $|A|$ of a set $A \subseteq G$ which does not contain $k$ distinct elements in arithmetic progression.

 The famous theorem of Szemer\'edi essentially asserts that $r_k(\Z/N\Z) = o_k(N)$. It is known, in fact, that the estimate $r_k(G) = o_k(N)$ holds for all $G$.

 There have been many papers concerning the issue of finding quantitative bounds for $r_k(G)$. A result of Bourgain states that 
\[ r_3(G) \ll N(\log \log N/\log N)^{1/2}\]
for all $G$. In this paper we obtain a similar bound for $r_4(G)$ in the particular case $G = F^n$, where $F$ is a fixed finite field with $\mbox{char}(F) \neq 2,3$ (for example, $F = \mathbb{F}_5$). We prove that
\[ r_4(G) \ll_{F} N(\log N)^{-c}\]
for some absolute constant $c > 0$. In future papers we will treat general abelian groups $G$, eventually obtaining a comparable result for arbitrary $G$. 

\end{abstract}

\maketitle

\section{Introduction}

 Let $G$ be a finite abelian group with cardinality $N$, written additively. Let $k \geq 3$ be an integer, and suppose that $(N,(k-1)!) = 1$ (or equivalently, that every non-zero element of $G$ has order at least $k$). We define $r_k(G)$ to be the largest cardinality $|A|$ of a set $A \subseteq G$ which does not contain an arithmetic progression $(x, x+ d, \dots x + (k-1) d)$ with $d \neq 0$ (such progressions will be referred to as \emph{proper}). 

 A deep and famous theorem of Szemer\'edi \cite{szemeredi} asserts that any set of integers with positive upper density contains arbitrarily long arithmetic progressions. This is easily seen to be equivalent to the assertion that 
\begin{equation}\label{eq0.00} r_k(\Z/N\Z) = o_k(N).\end{equation}
Here $o_k(N)$ denotes a quantity which when divided by $N$, goes to zero as $N \to \infty$ for each fixed $k$.
It is known, in fact, that $r_k(G) = o_k(N)$ for all $G$; this may be proved by combining Szemer\'edi's theorem with the density Hales-Jewett theorem \cite{fk}, and also follows from any of the recent hypergraph regularity results (see \cite{gowers-hyper}, \cite{rodlsurvey} and subsequent papers by the same authors and \cite{tao:hyper}).
When $k = 3$, the assertion \eqref{eq0.00} was proved earlier by Roth \cite{roth}, who in fact obtained the quantitative bound
\begin{equation}\label{eq0.01} r_3(\Z/N\Z) \ll N/\log \log N.\end{equation}
As usual we write $X \ll Y$ if we have the bound $X \leq CY$ for some absolute constant $C > 0$; if this constant $C$ depends on some additional parameters then we will denote this by subscripting the $\ll$ notation appropriately.
Roth's bound was improved to 
\begin{equation}\label{eq0.02} r_3(\Z/N\Z) \ll N(\log N)^{-c},\end{equation}
for some absolute constant $c > 0$, independently by Heath-Brown \cite{heath} and Szemer\'edi \cite{szemeredi-4}. This bound was then further improved to 
\begin{equation}\label{eq0.03} r_3(\Z/N\Z) \ll N(\log \log N/\log N)^{1/2}\end{equation} by Bourgain \cite{Bou}. This is the best bound currently known. It should be compared with the famous conjecture of Erd\H{o}s and Tur\'an \cite{erdos}, which asserts that if $A \subseteq \N$ is a set of integers with $\sum_{a \in A} a^{-1} = \infty$, then $A$ contains progressions of length $k$ for all $k$. This statement is unknown even when $k = 3$, in which case it is more-or-less equivalent to proving that $r_3(\Z/N\Z) \ll_\eps N/(\log N)^{1 + \eps}$ for all $\eps > 0$.

 Finding quantitative bounds for $r_4(\Z/N\Z)$ proved to be much more difficult. Many of the known proofs that $r_4(\Z/N\Z) = o(N)$, such as \cite{furst,gowers-hyper-4,rodl3,roth-4,szemeredi-4,tao:ergodic,taovu-book}, give very weak bounds or no explicit bounds at all. It was a great advance when Gowers \cite{gowers-4-aps,gowers-long-aps} proved that 
\[ r_4(\Z/N\Z) \ll N(\log \log N)^{-c}\]
for some absolute constant $c > 0$. This is the best bound currently known. Our goal in this paper and the next two in the series is to bring our knowledge concerning $r_4(G)$ into line with that concerning $r_3$.

 The arguments of Roth, Heath-Brown, Szemer\'edi and Bourgain can all (with varying degrees of effort) be adapted to give bounds for $r_3(G)$ of the same strength as \eqref{eq0.01}, \eqref{eq0.02} and \eqref{eq0.03} above for a general $G$. It was observed in \cite{mesh-ap} that Roth's argument is particularly simple when $G = \mathbb{F}_3^n$. In fact in this setting all four of the arguments of \cite{Bou, heath,roth,szemeredi-4} are essentially the same and give the bound
\[ r_3(\mathbb{F}_3^n) \ll N/\log N.\]
This idea of looking at finite field models for additive problems has proved very fruitful. The chief reason for its success is that arguments of linear algebra are available in the finite field setting, but not in general groups (see the survey \cite{green-fin-field} for more information).

 Our main theorem in this paper is
\begin{theorem}\label{main}
Let $F$ be a fixed finite field with $\mbox{\emph{char}}(F) \neq 2,3$. Let $G = F^n$, and write $N := |F|^n$. Then we have the bound
\[ r_4(G) \ll N(\log_{|F|} N)^{-c},\]
for some absolute constant $c > 0$ \textup{(}in fact one can take $c = 2^{-21}$\textup{)}.  The implied constant is absolute.
\end{theorem}
\noindent\textit{Remark.} One might perhaps keep in mind the example $F = \mathbb{F}_5$.

 This paper, like all previous papers obtaining quantitative bounds for $r_k(G)$, uses the \textit{density increment strategy}. See \cite{green-fin-field} for a general discussion of this strategy, and the book \cite{taovu-book} for proofs of \eqref{eq0.01}, \eqref{eq0.02} and \eqref{eq0.03}. Gowers \cite{gowers-4-aps, gowers-icm} obtained his bound using this strategy and some \textit{quadratic Fourier analysis}. Indeed, his bound may be described as the quadratic version of Roth's argument for $r_3$. In \cite{gt-inverseu3} we obtained the bound
\begin{equation}\label{eq0.04} r_4(G) \ll N(\log \log N)^{-c}\end{equation}
by an elaboration of the same method. The argument for $G = \mathbb{F}_5^n$, which is rather simpler than the general case, may be found in \S 7 of that paper and contains some of the ideas which will be important later on.

 The first step of that argument, and indeed the one of the main results of \cite{gt-inverseu3}, was an \textit{inverse theorem for the Gowers $U^3(\mathbb{F}_5^n)$ norm}. Combined with a so-called \textit{generalized von Neumann theorem,} this implies a certain very useful dichotomy. Let $A \subseteq \mathbb{F}_5^n$ be a set with density $\alpha$. Then \textit{either} $A$ contains roughly $\alpha^4 N^2$ progressions of length four (and hence at least one non-trivial progression) \textit{or} $A$ has density at least $\alpha + c(\alpha)$ on some set of the form $\{ x \in \mathbb{F}_5^n: q(x) = \lambda \}$, where $q : \mathbb{F}_5^n \rightarrow \mathbb{F}_5$ is a quadratic form, and $c(\alpha) > 0$ is an explicit positive quantity depending only on $\alpha$. 

 The next step is to \textit{linearize} the level set $\{q(x) = \lambda\}$, covering it by cosets of a subspace of dimension about $n/2$. $A$ must have density at least $\alpha + c(\alpha)/2$ on at least one of these, and this gives the basis for a density increment argument.

 Linearization is very costly, and is the chief reason that the bound in \eqref{eq0.04} contains an iterated logarithm. One way of avoiding linearization would be to work the whole argument on joint level sets (``quadratic submanifolds'')
\begin{equation}\label{quad-submanifold} \{ x : q_1(x) = \lambda_1, q_2(x) = \lambda_2, \dots , q_d(x)= \lambda_d\},\end{equation}
obtaining density increments on successive sets of this type (with $d$ increasing at each stage). Obtaining the relevant $U^3$ inverse and generalized von Neumann theorems turns out to be very troublesome, though it can be done; we hope to report further on this strategy in a future paper.

 In this paper we adopt a compromise approach, which may be thought of as the quadratic analogue of the Heath-Brown and Szemer\'edi bound for $r_3$. Very roughly, we prove that \textit{either} $A$ has roughly $\alpha^4 N^2$ four-term APs \textit{or} there are some quadratics $q_1,\dots,q_d$ such that $A$ has density at least $\alpha + c'(\alpha)$ on a quadratic submanifold such as \eqref{quad-submanifold}. Here $c'(\alpha)$ is to be thought of as rather larger than $c(\alpha)$. Only now do we linearize, covering the quadratic submanifold by cosets of some subspace of (it turns out) dimension about $n/(d+1)$. Note that if we linearized the quadratics one at a time we would pass to a subspace of dimension $n/2^d$. The relative efficiency of linearizing several quadratics at a time, together with the larger density increment $c'(\alpha)$, is what leads to the improved bound of Theorem \ref{main}.

 We are indebted to Timothy Gowers for inspiring this project, which was in fact the starting point for our collaboration, preceding (and leading to) such results as \cite{green-tao-primes,gt-primes-4}.  

\section{General notation}\label{notation-sec}

 Let $A$ be a finite non-empty set and let $f: A \to \C$ be a function. We use the 
averaging notation \[\E_A(f) = \E_{x \in A} f(x) := \frac{1}{|A|} \sum_{x \in A} f(x).\]
More complex expressions such as $\E_{x \in A, y \in B} f(x,y)$ are similarly defined.  We also define the $L^p$ norms
$$ \|f\|_{L^p(A)} := \E_A(|f|^p)^{1/p}$$
for $1 \leq p < \infty$, with the usual convention $\|f\|_{L^\infty(A)} := \sup_{x \in A} |f(x)|$.  We also use the complex inner product
$$ \langle f, g \rangle_{L^2(A)} := \E_A( f \overline{g} ).$$
We say that $f$ is $1$-\emph{bounded} if $\|f\|_{L^\infty(A)} \leq 1$.

 If $A,B$ are finite sets with $B$ non-empty, we use $\P_B(A) := \frac{|A \cap B|}{|B|}$ to denote the density of $A$ in $B$.
If $A \subset W$ are finite sets, we use $1_A: W \to \R$ to denote the indicator function of $A$, thus $1_A(x) = 1$ when $x \in A$ 
and $1_A(x) = 0$ otherwise. We also write $1_{x \in A}$ for $1_A(x)$.  Thus for instance $\P_B(A) = \E_B(1_A)$ for all non-empty $B \subseteq W$.

\section{Affine geometry}\label{affine-sec}

 Observe that to prove Theorem \ref{main} it suffices to do so in the special case when $F$ has prime order, since a vector space over
a general finite field can also be interpreted as a vector space of equal or greater dimension over a field of prime order.

\textsc{important convention.} Henceforth $F$ will be a fixed finite field of prime order at least $5$. Without this assumption some of our later arguments (Lemma \ref{cwt}, for example) do not work properly.

 Theorem \ref{main} is stated in terms of a vector space over $F$.  It is convenient to have
an \emph{affine} perspective, so that our definitions are insensitive to the choice of origin and thus enjoy
a translation invariance.  In this section we recall some of the basic features of affine linear algebra.  The notation here
may appear somewhat excessive, but we present the material in this manner for pedagogical reasons, as we shall shortly be
developing quadratic analogues of many of the concepts in this section, using the same type of notation.

\begin{definition}[Affine spaces]  Let $G$ be a linear vector space over $F$.  We define an \emph{affine subspace} $W$ of $G$ to be
a translate of a linear subspace $\dot W$ of $G$ by some arbitrary coset representative $y \in W$.  We refer to $\dot W = W-W$ as 
the \emph{homogenization} of $W$ and $W$ as a \emph{coset} of $\dot W$; note that if $x \in W$ and $h \in \dot W$ then $x+h \in W$.  
Two affine spaces are said to be \emph{parallel} if they have the same homogenization.  We define the \emph{dimension} $\dim(W)$ of
an affine space to be the dimension of its homogenization, and if $W$ is an affine subspace of another affine space $W'$ we refer
to the quantity $\dim(W') - \dim(W)$ as the \emph{codimension} of $W$ in $W'$.
\end{definition}

\begin{definition}[Linear phase function]\label{lphase-def}  Let $W$ be an affine space.  An \emph{\textup{(}affine-\textup{)} linear phase function} on $W$ is any map $\phi: W \to \R/\Z$
to the circle group $\R/\Z$ (which we view additively) such that
$$ \phi(x+h_1+h_2) - \phi(x+h_1) - \phi(x+h_2) + \phi(x) = 0$$
for all $x \in W$ and $h_1,h_2 \in \dot W$.  
\end{definition}

 For any finite additive group $G$, define the \emph{Pontryagin dual} $\widehat G$ of $G$ to be the space of group homomorphisms $\xi: x \mapsto \xi \cdot x$
from $G$ to $\R/\Z$.  Given a linear phase function $\phi$ on $W$, we can define its \emph{gradient vector} $\nabla \phi \in \widehat{\dot W}$
by requiring the Taylor expansion
$$ \phi(x+h) = \phi(x) + \nabla \phi \cdot h$$
for all $x \in W$ and $h \in \dot W$; it is easy to verify that $\nabla \phi$ is well defined.  Also observe that if $\phi$ is a linear phase function
on $W$, then $\phi$ takes at most $|F|$ values, and the level sets of $\phi$ are parallel affine spaces of codimension at most $1$ in $W$.

 Every linear phase $\phi$ on $W$ defines an \emph{affine character} $e(\phi): W \to \C$, where $e: \R/\Z \to \C$ is the standard
homomorphism $e(x) := e^{2\pi i x}$.  These characters could be used to develop an ``affine-linear Fourier analysis''.  However it will be more
convenient to use traditional (non-affine) linear Fourier analysis.  Namely, if $G$ is a linear vector space and $f: G \to \C$ is a function, we
define the Fourier transform $\widehat f: \widehat G \to \C$ by the formula
$$ \widehat f(\xi) := \E_{x \in G} f(x) e(-\xi \cdot x).$$
Of course we have the Fourier inversion formula
$$ f(x) = \sum_{\xi \in \widehat G} \widehat f(\xi) e(\xi \cdot x)$$
and the Plancherel identity
$$ \|f\|^2_{L^2(G)} := \E_G |f|^2 = \sum_{\xi \in G} |\widehat f(\xi)|^2.$$

\section{The form $\Lambda$ and the $U^3(W)$ norm}\label{lambda-sec}

 Let us say that four affine spaces $W_0, W_1, W_2, W_3$ in a common ambient space $W'$ are \emph{in arithmetic progression} if
they are parallel with common homogeneous space $\dot W$, and if they form an arithmetic progression in the quotient space
$W'/\dot W$, or equivalently if there exist $x \in W'$ and $h \in \dot W'$ such that $W_j = x + jh + \dot W$ for all $j = 0,1,2,3$.
In particular for any affine space $W$, the quadruple $W,W,W,W$ is in arithmetic progression.

 In the Fourier-analytic or ergodic approaches to counting progressions of length $4$, 
a fundamental role is played by the quadrilinear form $\Lambda_{W_0,W_1,W_2,W_3}(f_0,f_1,f_2,f_3)$, defined for four affine spaces
$W_0,W_1,W_2,W_3$ in arithmetic progression together with functions $f_j: W_j \to \C$ by
\[ \Lambda_{W_0,W_1,W_2,W_3}(f_0,f_1,f_2,f_3) := \E_{x \in W_0, h \in W_1-W_0} f_0(x) f_1(x+h) f_2(x+2h) f_3(x+3h).\]
We shall abbreviate $\Lambda_{W_0,W_1,W_2,W_3}$ as $\Lambda$ when the spaces $W_0,W_1,W_2,W_3$ are clear from
context (in particular, the spaces $W_0,W_1,W_2,W_3$ will usually be equal).  This quantity is clearly related to the number of arithmetic progressions of length $4$ in a set $A \subseteq W$.  In particular,
if $A$ has no proper progressions of length four, then it is easy to see that
\begin{equation}\label{lambda4}
 \Lambda_{W,W,W,W}(1_A, 1_A, 1_A, 1_A) = \P_W(A) / |W|.
\end{equation}

 We can now describe our ``density increment'' step, which (as we shall shortly see) easily implies Theorem \ref{main} upon iteration. 

\begin{theorem}[Anomalous number of AP4s implies density increment]\label{dens-inc}
Let $W$ be an affine space, and let $f: W \to \R$ be a 1-bounded non-negative
function \textup{(}thus $0 \leq f(x) \leq 1$ for all $x \in W$\textup{)}.  Set $\delta := \E_W(f)$.  Suppose that
\begin{equation}\label{law}
|\Lambda_{W,W,W,W}(f,f,f,f) - \Lambda_{W,W,W,W}(\delta, \delta, \delta, \delta)| > \delta^4 / 2.
\end{equation}
Let $C_1 := 2^{20}$.  Then at least one of the following two statements hold:
\begin{itemize}
\item \textup{(}medium-sized density increment on large space\textup{)}  There exists an affine subspace $W'$ of $W$ with dimension satisfying 
\begin{equation}\label{med-1} \dim(W') \geq \dim(W) - (2/\delta)^{C_1} \end{equation}
such that we have the density increment
\begin{equation}\label{med-2} \E_{W'}(f) \geq \delta + 2^{-44}\delta^{16}.\end{equation}
\item \textup{(}Large density increment on medium-sized space\textup{)} There exists an integer $K$, $1 \leq K \leq 2^{33}\delta^{-12}$, and an affine subspace $W'$ of $W$ with dimension satisfying
\begin{equation}\label{lar-1} \dim(W') \geq \frac{1}{K+1} \dim(W) - (2/\delta)^{C_1} \end{equation}
such that we have the density increment
\begin{equation}\label{lar-2} \E_{W'}(f) \geq \delta(1  +  2^{-15}K^{1/3}).\end{equation}
\end{itemize}
\end{theorem}

\noindent\textit{Proof of Theorem \ref{main} assuming Theorem \ref{dens-inc}.} It suffices to prove the following claim, which is more-or-less equivalent to Theorem \ref{main}.

\textit{Claim.} Let $A$ be a subset of an affine space $W$ with $\P_W(A) = \delta$.  If
\[\dim(W) \geq (2/\delta)^{C_2},\]
where $C_2 := 2^{21}$, then $A$ contains a proper arithmetic progression of length four.

 We prove this by induction on $\dim(W)$. This induction may alternatively be viewed as an iterated application of Theorem \ref{dens-inc}.  We may assume that $\dim(W) \geq 2^{C_2}$ since the claim is vacuous otherwise. This provides a start for the induction.
Let $f := 1_A$, so that $\delta = \E_W(f)$. Supposing for a contradiction that $A$ does not contain 
proper arithmetic progressions of length four, we see from \eqref{lambda4} that
\[ \Lambda_{W,W,W,W}(f,f,f,f) = \delta / |F|^{\dim(W)} \leq \delta/5^{(2/\delta)^{C_2}} < \delta^4/2,\] and thus \eqref{law} holds.  Applying Theorem \ref{dens-inc}, we conclude that either 
\begin{enumerate}
\item[(i)] there is a medium-sized density increment on a large space, meaning that \eqref{med-1} and \eqref{med-2} hold, or
\item[(ii)] there is a large density increment on a medium-sized space, meaning that \eqref{lar-1} and \eqref{lar-2} both hold for some parameter $K$, $1 \leq K \leq 2^{33}\delta^{-12}$.
\end{enumerate}

 Suppose that (i) holds.  From \eqref{med-2} and the fact that $\E_W f = \delta$ we see that $\dim(W') < \dim(W)$. Applying the induction hypothesis, we see that it suffices to check that  
\[ (2/\delta)^{C_2} - (2/\delta)^{C_1} \geq (2/(\delta + 2^{-44} \delta^{16}))^{C_2},\]
or in other words that
\[ (1 + 2^{-44} \delta^{15})^{-C_2} \leq 1 - (\delta/2)^{C_2-C_1}.\]
Using the inequality $(1 + x)^{-a} \leq 1 - ax + \frac{1}{2} a(a+1)x^2$, valid for $x,a \geq 0$, this is easy to verify since $C_2-C_1$ and $C_2$ are so large.  

 Suppose alternatively that (ii) holds.  Applying the induction hypothesis once more, it suffices to show that
\[ \frac{1}{K+1} (2/\delta)^{C_2} - (2/\delta)^{C_1} \geq (2/\delta(1 + 2^{-14} K^{1/3} ))^{C_2},\]
or in other words that
\[ \frac{1}{(1 + 2^{-14} K^{1/3})^{C_2}} \leq \frac{1}{K+1} - \bigg(\frac{\delta}{2}\bigg)^{C_2-C_1}.\]
which is again easy to verify since $C_2, C_2-C_1$ are so large.\endproof

 Much of the material in subsequent sections revolves around the estimation of $\Lambda$ in various ways.  Let us present just two simple
estimates here.  If $f_0,f_1,f_2,f_3$ are 1-bounded functions on an arithmetic progression of affine spaces $W_0,W_1,W_2,W_3$, then we have
\begin{equation}\label{l1-gvn}
|\Lambda_{W_0,W_1,W_2,W_3}(f_0,f_1,f_2,f_3)| \leq \min_{0 \leq j \leq 3} \|f_j\|_{L^1(W_j)}.
\end{equation}
This follows immediately from applying the change of variables $x \mapsto x - jh$ followed by the triangle inequality.  It leads to
an easy consequence:

\begin{lemma}[$L^1$ controls $\Lambda$]\label{l1-lambda4} 
Let $f_0,f_1,f_2,f_3,g_0, g_1,g_2,g_3: W \to \C$ be functions on an affine space $W$ which are all uniformly bounded by some $\alpha > 0$.
Then we have
\[ |\Lambda_{W,W,W,W}(f_0,f_1,f_2,f_3) - \Lambda_{W,W,W,W}(g_0,g_1,g_2,g_3)| \leq 4 \alpha^3 \sup_{0 \leq i\leq 3} \|f_i-g_i\|_{L^1(W)}.\]
\end{lemma}

\begin{proof} By dividing $f_i$ and $g_i$ by $\alpha$ we may normalize so that $\alpha=1$.  We abbreviate $\Lambda_{W,W,W,W}$ as $\Lambda$.
The claim then follows from the telescoping identity
\begin{equation}\label{telescoping}
\begin{split}
\Lambda(f_0,f_1,f_2,f_3) - \Lambda(g_0,g_1,g_2,g_3) &= \Lambda(f_0-g_0,g_1,g_2,g_3) + \Lambda(f_0,f_1-g_1,g_2,g_3)\\
&\qquad 
+ \Lambda(f_0,f_1,f_2-g_2,g_3) + \Lambda(f_0,f_1,f_2,f_3-g_3)
\end{split}
\end{equation}
and \eqref{l1-gvn}.
\end{proof}

 The next lemma involves, for the first time in this paper, the \textit{Gowers $U^3$-norm} on $W$ (cf. \cite{gowers-4-aps,gowers-long-aps,green-tao-primes}). If $f : W \rightarrow \C$ is a function, recall that
\begin{eqnarray*} \Vert f \Vert_{U^3(W)}^8 & := & \E_{x \in W; h_1,h_2,h_3 \in \dot W}(f(x)\overline{f(x+h_1)f(x+h_2)f(x+h_3)}f(x+h_1+h_2) \times \\
& & \qquad\qquad\qquad\qquad \times f(x+h_2+h_3)f(x+h_1+h_3) \overline{f(x+h_1+h_2+h_3)}).\end{eqnarray*}
The Gowers $U^3$-norm measures the extent to which $f$ behaves ``quadratically''. Note for example that if $f(x) = \omega^{\phi(x)}$, where $\phi : F^n \rightarrow F$ is a quadratic form and $\omega = e(2\pi i/|F|)$, then $\Vert f \Vert_{U^3} = 1$, the largest possible $U^3$ norm of a 1-bounded function. The $U^3$-norm also controls progressions of length four in a sense to be made precise in Lemma \ref{gvn-lem} below. There are also Gowers $U^d$-norms for $d = 2,3,\dots$, with the $U^d$-norm controlling progressions of length $d+1$.
Properties of the Gowers norms may be found in \cite{gowers-4-aps,gowers-long-aps,green-tao-primes}: the paper \cite{gt-inverseu3} and book \cite{taovu-book} provide a comprehensive discussion of the $U^3$-norm. For a discussion in an ergodic theory context see \cite{host-kra}. Previous papers only consider the Gowers norms on abelian groups, but the generalization to affine spaces is a triviality.

\begin{lemma}[Generalized von Neumann]\label{gvn-lem} Suppose that $f_0,f_1,f_2,f_3 : W \rightarrow \C$ are 1-bounded functions. Then we have
\[ |\Lambda_{W,W,W,W}(f_0,f_1,f_2,f_3)| \leq \min_{0 \leq i \leq 3} \|f_i\|_{U^3(W)}.\]
In fact more generally if $W_0,W_1,W_2,W_3$ are in arithmetic progression and if $f_i : W_i \to \C$ are functions then we have
\begin{equation}\boxeq |\Lambda_{W_0,W_1,W_2,W_3}(f_0,f_1,f_2,f_3)| \leq \min_{0 \leq i \leq 3} \|f_i\|_{U^3(W_i)}.\end{equation}
\end{lemma}
\textit{Remarks.} The first statement is \cite[Proposition 1.7]{gt-inverseu3}, and is proved in \S 4 of that paper using three applications of the Cauchy-Schwarz inequality. Versions of this inequality appear in several earlier works also, such as \cite{gowers-4-aps}.  The second statement may be proved in the same way, with trivial notational changes.

 Using the telescoping identity \eqref{telescoping}, we conclude the following variant of Lemma \ref{l1-lambda4}:

\begin{lemma}[$U^3$ controls $\Lambda$]\label{u3-lambda4} 
Let $f, g: W \to \C$ be 1-bounded functions on an affine space $W$.  Then we have
$$ |\Lambda_{W,W,W,W}(f,f,f,f) - \Lambda_{W,W,W,W}(g,g,g,g)| \leq 4 \|f-g\|_{U^3(W)}.$$
\end{lemma}

 Lemmas \ref{l1-lambda4} and \ref{u3-lambda4} show that errors with small $L^1$ or $U^3$ norm are negligible for the purposes of counting
progressions of length $4$.  We understand the $L^1$ norm very well, but the $U^3$ norm is far more mysterious. For us, a key tool in its study will be the \textit{inverse theorem} of \cite{gt-inverseu3}, providing a description of those $f$ for which $\Vert f \Vert_{U^3}$ is large.

\section{The inverse $U^3(W)$ theorem}\label{inverse-sec}

 We begin with some notation.

\begin{definition}[Quadratic phase function]  Let $W$ be an affine space.  An \emph{\textup{(}affine\textup{)} quadratic phase function} on $W$ is any map $\phi: W \to \R/\Z$
such that
\begin{align*}
 \phi(x+h_1+h_2+h_3) &- \phi(x+h_1+h_2) - \phi(x+h_2+h_3) - \phi(x+h_1+h_3) \\
 &+ \phi(x+h_1) + \phi(x+h_2) + \phi(x+h_3) - \phi(x) = 0
\end{align*}
for all $x \in W$ and $h_1,h_2,h_3 \in \dot W$.  
\end{definition}

 Let us make the trivial remark that the translation of a quadratic phase function remains quadratic (even if one also translates the underlying space $W$), and the restriction of a quadratic phase function to an affine subspace remains quadratic.  Also, every linear phase function is automatically
quadratic. An example of a quadratic phase function on $G = F^n$ is $\phi(x) = \nu(Mx \cdot x + \xi \cdot x) + c$, where $M : F^n \rightarrow F^n$ is a self-adjoint linear transformation, $\xi \in \widehat{G}$, $c \in \R/\Z$ and $\nu : F \rightarrow \R/\Z$ is some fixed homomorphism of additive groups. In fact (as we shall see) every quadratic phase can be written explicitly in this way.

 Quadratic phases are 
closely tied to the $U^3(W)$ norm; indeed one can easily verify that a 1-bounded function $f : W \to \C$ has $U^3(W)$ norm bounded by $1$, with equality holding
if and only if $f = e(\phi)$ for some quadratic phase $\phi : W \rightarrow \R/\Z$.  A more quantitative version of this fact is as follows.

\begin{theorem}[Inverse theorem for $U^3(W)$]\label{invert-u3}\cite[Theorem 2.3]{gt-inverseu3}  
Let $f: W \to \C$ be a 1-bounded function on an affine space $W$
such that $\|f\|_{U^3(W)} \geq \eta$ for some $0 < \eta \leq 1$.  Let $C_0 := 2^{16}$.
Then there exists a linear subspace $\dot W'$ of $\dot W$ of codimension at most $(2/\eta)^{C_0}$	
such that for each coset $W'$ of $\dot W'$, there exists a quadratic phase function $\phi_{W'}: W' \to \R/\Z$ such that
\begin{equation}\label{ygw}
 \E_{W' \in W/\dot W'} |\E_{x \in W'} f(x) e(- \phi_{W'}(x) )| \geq (\eta/2)^{C_0}.
 \end{equation}
\end{theorem}

\textit{Remarks.} The result in \cite{gt-inverseu3} is phrased for vector spaces over $F$ rather than affine spaces but the extension to the
affine case is a triviality.  Also the averaging in \cite{gt-inverseu3} is over coset representatives rather than actual cosets but the two
are related by an easy application of the pigeonhole principle.
There is a corresponding theorem in arbitrary finite groups $G$ but it is somewhat more complicated in that the subspace $W$ needs
to be replaced by a Bohr set: see \cite{gt-inverseu3}.  An analogue of this result also holds in the characteristic $2$ case
(Samorodnitsky, private communication) but we will not need it here.

 In \cite{gt-inverseu3} it was conjectured that one could in fact take $\dot W'=\dot W$ (possibly at the cost of deteriorating the constant 
$C_0 = 2^{16}$), in which
case this inverse theorem takes a particularly simple form\footnote{One can use a simple Fourier averaging, combined with a certain ``quadratic extension theorem'' to achieve a version of this, but with bounds that deteriorate exponentially in $\eta$, see \cite{gt-inverseu3}.  If we wish to prove Theorem \ref{main} we cannot afford such exponential losses and so will not use this fact.}: if the $U^3(G)$ norm of $f$ is large, then $f$ has large correlation with $e(\phi)$ for some quadratic phase $\phi : W \to \R/\Z$.  This would simplify the arguments in this paper somewhat, though in practice it is relatively inexpensive to pass from $W$ to the slightly smaller space $W'$ as necessary. 

 As mentioned in the introduction, one can use this theorem (together with Corollary \ref{linquad} below)
to run a density increment argument. This yields
a weak version of Theorem \ref{dens-inc}, giving a
bound of the form $r_4(F^n) \ll_F N(\log \log N)^{-c}$ for some $c > 0$: see \cite[\S 7]{gt-inverseu3} for the details.  We will however take a more efficient
route involving the \emph{energy increment argument} from \cite{green-tao-primes}, motivated both by considerations from ergodic
theory (notably Furstenberg's ergodic proof \cite{furst,FKO} of Szemer\'edi's theorem) and from regularity lemmas from graph theory
(in particular Szemer\'edi's regularity lemma \cite{szemeredi}, as well as an arithmetic analogue of the first author \cite{green1}).

\section{Linear and quadratic factors, and a quadratic Koopman-von Neumann theorem}\label{kvn-sec}

 To convert the inverse theorem to a quadratic structure theorem (or quadratic Koopman-von Neumann theorem)
we need some concepts from ergodic theory. Note, however, that in our context all of these constructions are purely finitary.

\begin{definition}[Factors]  Let $W$ be an arbitrary finite non-empty set (typically $W$ will be an affine space).
Define a \emph{factor} (or \emph{$\sigma$-algebra})
in $W$ to be a collection $\B$ of subsets of $W$ which are closed under union, intersection, and complement, and which contains $\emptyset$ and $W$. Define an \emph{atom} of $\B$ to be a minimal non-empty element of $\B$; these partition $W$, and indeed in this finitary setting the set of factors
can be placed in one-to-one correspondence with the set of partitions of $W$.  A function $f: W \to X$ with arbitrary finite range $X$
is said to be \emph{measurable with respect to $\B$}, or
$\B$-measurable for short, if all its level sets lie in $\B$.  We let $\B_f$ denote the factor generated by $f$, thus the atoms of $\B_f$ are precisely the non-empty level sets of $f$.  We say that one factor $\B'$ \emph{extends} another $\B$ if $\B \subset \B'$; we also say that $\B$ is a \emph{factor} of $\B'$ in this case.
If $W'$ is any non-empty subset of $W$ and $\B$ is a factor in $W$, we define the \emph{restriction} $\B|_{W'}$ of $\B$ to $W'$ to be the
factor in $W'$ formed by intersecting all the sets in $\B$ with $W'$.  Note that this is a subset of $\B$ if $W' \in \B$.
If $\B, \B'$ are factors in $W$ we let $\B \vee \B'$ be the smallest common extension (thus the atoms of $\B \vee \B'$ are the intersections of atoms of $\B$ and atoms of $\B'$).  If $f: W \to \C$, we let $\E(f|\B): W \to \C$ denote the conditional expectation
$$ \E(f|\B)(x) := \E(f|\B(x)) \hbox{ for all } x \in W,$$
where $\B(x)$ is the unique atom in $\B$ that contains $x$.  Equivalently, $\E(f|\B)$ is the orthogonal projection (in the Hilbert space
$L^2(W)$) of $f$ to the space of $\B$-measurable functions.  
\end{definition}

 We will focus our attention on very structured factors, namely linear and quadratic factors, which turn out in the finite field setting to be
the only factors required to analyze progressions of length four.  

\begin{definition}[Linear factors]  Let $W$ be an affine space.  A \emph{linear factor of complexity at most $d$} is any factor $\B$ in $W$ 
of the form $\B = \B_{\phi_1} \vee \ldots \vee \B_{\phi_{d'}}$, where $0 \leq d' \leq d$ and $\phi_1, \ldots, \phi_d$ are linear phase functions on $W$.
\end{definition}

 Observe that if $\B$ is a linear factor of complexity at most $d$, then the atoms of $\B$ are parallel affine spaces of codimension at most $d$.
Also, if $\B$ and $\B'$ are linear factors of complexity at most $d$, $d'$, then $\B \vee \B'$ is also a linear factor, with complexity at most $d+d'$.

\begin{definition}[Quadratic factors]  Let $W$ be an affine space.  A \emph{pure quadratic factor of complexity at most $d$} is
any factor $\B$ in $W$ of the form $\B = \B_{\phi_1} \vee \ldots \vee \B_{\phi_{d'}}$, where $0 \leq d' \leq d$ and $\phi_1, \ldots, \phi_d$ are quadratic phase functions on $W$.  A \emph{quadratic factor of complexity at most $(d_1,d_2)$} is any pair $(\B_1,\B_2)$ of factors in $W$, where $\B_1$
is a linear factor of complexity at most $d_1$, and $\B_2$ is an extension of $\B_1$, whose restriction to any atom of $\B_1$ is a pure quadratic factor of complexity at most $d_2$.  We say that one quadratic factor $(\B'_1,\B'_2)$ is a \emph{quadratic extension} of
another $(\B_1,\B_2)$ if $\B_1 \subseteq \B'_1$ and $\B_2 \subseteq \B'_2$.
\end{definition}

\textit{Remark.} Note that a linear factor of complexity $d_1$ can have as many as $|F|^{d_1}$ atoms, and thus a quadratic factor of complexity $(d_1,d_2)$ can involve as many as $|F|^{d_1} d_2$ quadratics (on up to $|F|^{d_1}$ different domains), though it is quite likely that an improved version of the inverse theorem in \cite{gt-inverseu3} would reduce the number of quadratics involved here.  Some care must be taken to avoid 
this exponential dependence on the complexity from destroying the polynomial nature of many of the quantities in our arguments.  Fortunately, by 
working ``locally'' on linear atoms rather than globally on all of $W$ one can avoid any unpleasant factors of $|F|^{d_1}$ in our analysis.

 Observe that if $(\B_1,\B_2)$ and $(\B'_1,\B'_2)$ are quadratic factors of complexity at most $(d_1,d_2)$ and $(d'_1,d'_2)$ respectively,
then their common extension $(\B_1 \vee \B'_1, \B_2 \vee \B'_2)$ is a quadratic factor of complexity at most $(d_1+d'_1,d_2+d'_2)$; this is 
because the restriction of a pure quadratic factor to an affine subspace remains a pure quadratic factor of equal or lesser complexity.

 The inverse theorem, Theorem \ref{invert-u3}, can now be rephrased in terms of quadratic factors as follows.

\begin{theorem}[Inverse theorem for $U^3(W)$, again]\label{invert-u3-sigma}
Let $f: W \to \C$ be a 1-bounded function on an affine space $W$
such that $\|f\|_{U^3(W)} \geq \eta$ for some $\eta$, $0 < \eta \leq 1$.  Then there exists a quadratic factor $(\B_1,\B_2)$
in $W$ of complexity at most $((2/\eta)^{C_0}, 1)$ such that
$$ \| \E(f|\B_2) \|_{L^1(W)} \geq (\eta/2)^{C_0}.$$
\end{theorem}

 This has the following consequence.

\begin{corollary}[Lack of relative uniformity implies energy increment]\label{energy-increment}  
Let $(\B_1,\B_2)$ be a quadratic factor of complexity at most $(d_1,d_2)$
in an affine space $W$, and let $f: W \to \R^+$ be a 1-bounded non-negative function such that
$\| f - \E(f|\B_2)\|_{U^3(W)} \geq \eta$ for some $\eta$, $0 < \eta \leq 1$.  There exists a quadratic extension $(\B'_1,\B'_2)$ of $(\B_1,\B_2)$
of complexity at most $(d_1+ (2/\eta)^{C_0}, d_2+1)$ such that we have the energy increment
$$ \| \E(f|\B'_2) \|_{L^2(W)}^2 \geq \| \E(f|\B_2) \|_{L^2(W)}^2 + (\eta/2)^{2C_0}.$$
\end{corollary}

\begin{proof} Applying Theorem \ref{invert-u3-sigma} to the 1-bounded\footnote{Here we are using the hypothesis that $f$ is non-negative and bounded by $1$ to ensure that $\Vert f - \E(f | \B_2) \Vert_{\infty} \leq 1$. More generally we would replace the $1$ on the right-hand side by a $2$, which has a negligible impact on the argument.} function $f - \E(f|\B_2)$, we can find a quadratic factor $(\widetilde \B_1, \widetilde \B_2)$ of complexity at most $((2/\eta)^{C_0},1)$ such that
$$ \| \E( f - \E(f|\B_2) | \widetilde \B_2) \|_{L^1(W)} \geq (\eta/2)^{C_0}.$$
If we then let $\B'_1 := \B_1 \vee \widetilde \B_1$ and $\B'_2 := \B_2 \vee \widetilde \B_2$, then by Pythagoras' theorem, the inclusions
$\B_2, \widetilde \B_2 \subseteq \B'_2$, and Cauchy-Schwarz, we have
\begin{align*}
\| \E(f|\B'_2) \|_{L^2(W)}^2 - \| \E(f|\B_2) \|_{L^2(W)}^2
&= \| \E(f|\B'_2) - \E(f|\B_2) \|_{L^2(W)}^2 \\
&= \| \E( f - \E(f|\B_2) | \B'_2 ) \|_{L^2(W)}^2 \\
&\geq \| \E( f - \E(f|\B_2) | \widetilde \B_2 ) \|_{L^2(W)}^2 \\
&\geq \| \E( f - \E(f|\B_2) | \widetilde \B_2 ) \|_{L^1(W)}^2 \\
&\geq (\eta/2)^{2C_0}.
\end{align*}
The claim follows.
\end{proof}

 We can employ this corollary repeatedly. This ``energy increment argument'' allows us to deduce one of our most important tools, a \textit{\textup{(}Quadratic\textup{)} Koopman-von Neumann theorem}\footnote{This theorem decomposes a function $f$ orthogonally into a ``quadratically almost periodic'' component $\E(f|\B_2)$ and a ``quadratically mixing'' component $f - \E(f|\B_2)$.  This can be compared with the ordinary (linear) Koopman-von Neumann theorem in (infinitary) ergodic theory, which splits a function $f$ orthogonally into an almost periodic component and a weakly mixing component.  See also the analysis of characteristic factors for the Gowers norms and for multiple recurrence in \cite{bhk,host-kra,ziegler}.}.

\begin{theorem}[Quadratic Koopman-von Neumann theorem]\label{kvn}  
Let $f: W \to \C$ be a 1-bounded non-negative function on an affine space $W$, and let $\eta > 0$.  Then there exists a
quadratic factor $(\B_1,\B_2)$ in $W$ of complexity at most $((2/\eta)^{3C_0}, (2/\eta)^{2C_0})$
such that
\begin{equation}\label{efb-b1}
\|f-\E(f|\B_2)\|_{U^3(W)} \leq \eta.
\end{equation}
\end{theorem}

\begin{proof} Start with $(\B_1,\B_2) = (\{ \emptyset, W\}, \{\emptyset,W\})$, which is a quadratic factor of complexity $(0,0)$.  If
\eqref{efb-b1} holds then we are done. Otherwise, we may apply Corollary \ref{energy-increment} to extend $(\B_1,\B_2)$ to a quadratic factor
with complexity incremented by at most $((2/\eta)^{C_0}, 1)$ and the energy 
$\| \E(f|\B_2)\|_{L^2(W)}^2$ incremented by at least $(\eta/2)^{2C_0}$. On the other hand,
since $f$ is 1-bounded, the energy $\| \E(f|\B_2)\|_{L^2(W)}^2$ is positive and at most $1$.  Thus we cannot iterate the above 
procedure more than $(2/\eta)^{2C_0}$ times before terminating.  The claim follows.
\end{proof}

 Applying this result for $\eta := \delta^4/16$ and then using Lemma \ref{u3-lambda4}, we conclude

\begin{corollary}[Too few AP4s on a quadratic factor]\label{dens-cor}
Let $W$ be an affine space, and let $f: W \to \R$ be a 1-bounded non-negative
function.  Set $\delta := \E_W(f)$.  Suppose that
$$ |\Lambda_{W,W,W,W}(f,f,f,f) - \Lambda_{W,W,W,W}(\delta, \delta, \delta, \delta)| > \delta^4 / 2.$$
Then there exists a quadratic factor $(\B_1,\B_2)$ in $W$ of complexity at most $(d_1,d_2)$, where $d_1  := (32/\delta^4)^{3C_0}$ and $d_2 :=  (32/\delta^4)^{2C_0}$,
such that the function $g := \E(f|\B_2)$ obeys
$$ |\Lambda_{W,W,W,W}(g,g,g,g) - \Lambda_{W,W,W,W}(\delta, \delta, \delta, \delta)| > \delta^4 / 4.$$
\end{corollary}

 The factor $(\B_1,\B_2)$ is closely related to the ergodic theory concept of a \emph{characteristic factor} for the problem of obtaining
$4$-term recurrence for $f$.  The corollary thus replaces the study of $f$ (which could essentially be an arbitrary function) by
a much lower complexity object $g$, which in principle can be described explicitly using a bounded number of linear and quadratic
phase functions (cf. \cite{bhk}).  It will be the first component of the proof of Theorem \ref{dens-inc}.

 Of course, it remains to understand $g$, and more precisely to count progressions of length four in $g$.  This will be second component of
the proof of Theorem \ref{dens-inc}, and will require a certain amount of (fairly standard) analysis of the geometry, algebra, and Fourier-analytic
structure of quadratic phase functions and their associated level sets.  This will be the objective of the remaining sections of the paper.

\section{Affine quadratic geometry}\label{geometry-sec}

 In order to analyze $g$, we first must understand the geometry of quadratic phase functions.
From an algebraic perspective, at least, the structure of quadratic phase functions can be easily understood.
Given a quadratic phase function $\phi$ on an affine space $W$, define the \emph{gradient} $\nabla \phi: W \to \widehat{\dot W}$ and the 
\emph{Hessian} $\nabla^2 \phi: \dot W \to \widehat{\dot W}$ by requiring the Taylor expansion
\begin{equation}\label{taylor}
 \phi(x+h) = \phi(x) + \nabla \phi(x) \cdot h + \frac{1}{2} \nabla^2 \phi h \cdot h
 \end{equation}
for all $x \in W$ and $h \in \dot W$.  Indeed we have the explicit formulae
$$ \nabla^2 \phi h_1 \cdot h_2 = \phi(x+h_1+h_2) - \phi(x+h_1) - \phi(x+h_2) + \phi(x)$$
(note the right-hand side is in fact independent of $x$) and
$$ \nabla \phi(x) \cdot h = \frac{1}{2}( \phi(x+h) - \phi(x-h) );$$
here we exploit the hypothesis that $|F|$ is odd in order to be able to divide by $2$.
Observe that $\nabla \phi$ is linear and that $\nabla^2 \phi$ is a self-adjoint linear transformation from $\dot W$ to $\widehat{\dot W}$, which vanishes if and only if
$\phi$ is linear (in which case $\nabla \phi(x) = \nabla \phi$ is constant); this combined with \eqref{taylor}
shows that quadratic phase functions do indeed have the form $\phi(x) = \nu(Mx \cdot x + \xi \cdot x) + c$ as claimed earlier. Note that if one chooses a basis $(e_1,\dots,e_n)$ for $F^n$ then, in the associated coordinate system, $\phi$ takes the rather concrete form
\[ \phi(x_1,\dots,x_n) = \nu\big(\sum_{1 \leq i \leq j \leq n} M_{ij}x_ix_j + \sum_{i = 1}^n r_i x_i\big) + c.\] Furthermore, $\nabla^2 \phi$ has a null space
$\ker(\nabla^2 \phi)$, which is a linear subspace of $\dot W$; observe that $\phi$ becomes linear when restricted to any coset of this space.
The codimension of this null space will be referred to as the \emph{rank} $\rank(\phi)$ of $\phi$.  Intuitively,
this rank measures how close the quadratic phase function is to being linear.

 One can easily check (from \eqref{taylor}) that a quadratic phase function takes at most $|F|$ values; indeed, after shifting the 
phase by a constant phase in $\R/\Z$
we may assume that the quadratic phase function takes values in the discrete group $\T_F := \{ x \in \R/\Z: |F| x = 0 \}$.  Let us refer
to such phase functions as \emph{discretized}.  Note that the space of discretized quadratic phase functions is itself a finite-dimensional
vector space over $F$.

 We are now in a position to study the structure of quadratic factors, and in particular their atoms (which are nothing more than quadratic varieties).  
We begin with the classical result of Chevalley and Warning.

\begin{lemma}[Chevalley-Warning theorem]\label{cwt}  Let $W$ be an affine space, and let $\B$ be a pure quadratic factor of complexity strictly less than $\dim(W)/2$.  Then every atom of $W$ has size a multiple of $|F|$.  In particular, if an atom
contains one point $x_0 \in W$, then it must contain at least one further point in $W$.
\end{lemma}

\begin{proof}  We may translate $W$ to be a linear space, which we then identify with $F^n$.
Subtracting constant terms as appropriate, we may write any atom $A$ of $\B$ as
$$ A = \{ x \in F^n: \phi_1(x) = \ldots = \phi_d(x) = 0 \}$$
for some discretized quadratic phases $\phi_1,\ldots,\phi_d: W \to \T_F$ and some $d < n/2$.  
Identifying $\T_F$ with $F$ and writing $\phi_1,\ldots,\phi_d$ in coordinates,
it follows that
$$ A = \{ x \in F^n: Q_1(x) = \ldots = Q_d(x) = 0 \}$$
for some quadratic polynomials $Q_1,\ldots,Q_d: F^n \to F$.  Modulo $|F|$, we thus have
\begin{equation}\label{star-33} |A| \equiv \sum_{x \in F^n} \prod_{j=1}^d (1 - Q_j(x)^{|F|-1}) \md{|F|}.\end{equation}
The product has degree at most $2d(|F|-1)$. Since $d < n/2$, we see after writing $x = (x_1,\ldots,x_n)$ that none of the monomials
in this product are multiples of $x_1^{|F|-1} \ldots x_n^{|F|-1}$.  The right-hand side of \eqref{star-33} therefore vanishes, and we are done.
\end{proof}

 We now apply this lemma to obtain large linear spaces inside quadratic varieties.  We begin with a homogeneous statement.

\begin{lemma}[Quadratic forms have large null spaces]  Let $G$ be a vector space, and let 
$M_1,\ldots,M_d$ be self-adjoint linear transformations from $G$ to $\widehat G$.
Then there exists a linear subspace $\dot W$ of $G$ with \[\dim(\dot W) \geq \frac{1}{d+1} \dim(G) - \frac{2d}{d+1}\]
such that $M_j x \cdot y = 0$ for all $j$,  $1 \leq j \leq d$ and for all $x,y \in \dot W$.
\end{lemma}

\begin{proof}  Let $\dot W$ be a maximal linear subspace of $G$ which is null with respect to all the $M_j$ (i.e. $M_j x \cdot y = 0$ for all $j$, 
$1 \leq j \leq d$ and for all $x,y \in \dot W$).  Let $\dot W^\perp$ be the linear subspace 
$$\dot W^\perp := \{ x \in G: M_j x \cdot y = 0 \hbox{ for all } 1 \leq j \leq d \hbox{ and } y \in \dot W \}$$
thus $\dot W^\perp \supseteq\dot W$.  From linear algebra we also see that
$$ \dim(\dot W^\perp) \geq \dim(G) - d \dim(\dot W)$$
and thus after some algebra
$$ \dim(\dot W) \geq \frac{1}{d+1} \dim(G) - \frac{1}{d+1} \dim(\dot W^\perp/\dot W).$$
Observe that the quadratic forms $Q_j(x) := M_j x \cdot x$ are well-defined on $\dot W^\perp/\dot W$.  The zero locus $\{ x \in \dot W^\perp/\dot W: Q_j(x) = 0 \hbox{ for all $j$, } 1 \leq j \leq d \}$ consists only of the origin, since otherwise we could extend $\dot W$
by one additional dimension and contradict maximality.  In particular, the cardinality of this zero locus is not a multiple of $|F|$. 
Applying Lemma \ref{cwt} we conclude that $\dim(\dot W^\perp/\dot W) \leq 2d$, and the claim follows.
\end{proof}

\begin{corollary}[Linearization of quadratic factors]
Let $(\B_1,\B_2)$ be a quadratic factor on an affine space $W$ of complexity at most $(d_1,d_2)$.  Then each atom of $\B_2$ can be partitioned 
into disjoint affine spaces each of dimension at least
$$ \frac{\dim(W)}{d_2+1} - \frac{d_1+2d_2}{d_2+1} - d_2.$$
\end{corollary}

\begin{proof} By working on each atom of $\B_1$ separately, one sees that it suffices to verify this claim for pure quadratic factors.  Thus we may take $\B_1$ to be trivial and $\B_2 = \B_{\phi_1} \vee \ldots \vee \B_{\phi_{d_2}}$ for some quadratic
phase functions $\phi_1,\ldots,\phi_{d_2}$.
By the preceding lemma
we can find a linear subspace $\dot W'$ of $\dot W$ of dimension
\[ \dim(\dot W') \geq \frac{1}{d_2+1} \dim(W) - \frac{2d_2}{d_2+1} \]
which is null with respect to all of the $\nabla^2 \phi_j$.  In particular, this implies that $\phi_1,\ldots,\phi_{d_2}$ are 
linear on each of the cosets of $\dot W'$ in $W$.  Thus one can refine each such coset $W'$ further into affine spaces of dimension at most 
$\dim(\dot W') - d_2$, on which each of the $\phi_1,\ldots,\phi_{d_2}$ are constant.  These spaces form a partition of the atoms of $\B_2$, and the
claim follows.\end{proof}

 As a consequence of this, we see that a density increment on a quadratic factor implies a density increment on a subspace, albeit at the
expense of reducing the domain of the density increment substantially.

\begin{corollary}[Linearization of quadratic density increment]\label{linquad} Let $f: W \to \R$ be a real-valued function on an
affine space $W$, and let $(\B_1,\B_2)$ be a quadratic factor of complexity at most $(d_1,d_2)$. 
Let $A$ be an atom of $\B_2$.  Then there exists an affine subspace $W'$ of $W$ of dimension
$$ \dim(W') \geq \frac{1}{d_2+1} \dim(W) - \frac{d_1+2d_2}{d_2+1} - d_2$$
such that $\E_{W'}(f) \geq \E_A(f)$.
\end{corollary}

\begin{proof} From the preceding corollary, we can write $A$ as the disjoint union of affine spaces of
dimension at least $\frac{1}{d_2+1} \dim(W) - \frac{d_1+2d_2}{d_2+1} - d_2$.  The claim then follows from the pigeonhole
principle.
\end{proof}

 We now record a simple linear variant of this which will be used in \S \ref{mixing-sec}. 

\begin{lemma}\label{lininc} Let $f: W \to \R$ be a real-valued function on an
affine space $W$, and let $\B$ be a linear factor of complexity at most $d$.  Then there exists an affine subspace $W'$
of $W$ with $\dim(W') \geq \dim(W) - d$ such that $\E_{W'}(f) \geq \E_W(f) + \frac{1}{2}\| \E(f|\B) - \E_W(f) \|_{L^1(W)}$.
\end{lemma}

\begin{proof} Suppose that $g : W \rightarrow \R$ is any function with mean zero. Then we have $E(g + |g|) = \|g \|_1$, which implies that there is some $x \in W$ such that $g(x) + |g(x)| \geq \| g \|_1$. For such an $x$ we clearly have $g(x) \geq \frac{1}{2}\| g \|_1$. Applying this observation with $g := \E(f|\B)-\E_W(f)$, we see that there must exist $x \in W$ such that
\[ \E(f|\B)(x) - \E_W(f) \geq \frac{1}{2}\| \E(f|\B) - \E_W(f) \|_{L^1(W)}.\]  Letting $W'$ be the atom of $\B$ containing $x$, the claim follows.
\end{proof}

 At this point we can already conclude a cheap version of Theorem \ref{dens-inc}, the density increment result:

\begin{proposition}[Cheap density increment for 4APs]\label{cheap-inc}  
Let $W$ be an affine space, and let $f: W \to \R$ be a 1-bounded non-negative
function.  Set $\delta := \E_W(f)$.  Suppose that
$$ |\Lambda_{W,W,W,W}(f,f,f,f) - \Lambda_{W,W,W,W}(\delta, \delta, \delta, \delta)| > \delta^4 / 2.$$
Then there exists an affine subspace $W'$ of $W$ satisfying the dimension bound
\[
 \dim(W) \geq (\delta^4/32)^{2C_0+1} \dim(G) - (32/\delta^4)^{2C_0+1}
\]
and such that we have the density increment
\[
\E_{W'}(f) \geq \delta + \frac{\delta}{64}.
\]
\end{proposition}

\textit{Remark.} An iteration of the above proposition gives a bound of the form
\begin{equation}\label{cheap}
r_4(F^n) \ll Ne^{-c\sqrt{\log\log_{|F|}N}}
\end{equation}
for some absolute constant $c > 0$ (recall that $N := |F|^n$). We leave the verification of this to the reader, remarking that it is very similar to the deduction of Theorem \ref{main} from Theorem \ref{dens-inc} as given in \S \ref{lambda-sec}. The bound \eqref{cheap} is better than the previously best-known result, \eqref{eq0.04}, but substantially weaker than Theorem \ref{main}. It does enjoy the advantage of being easier to adapt to groups more general than $F^n$: for details see \cite{green-tao-szem2}.

\textit{Proof of Proposition \ref{cheap-inc}.}    
We apply Corollary \ref{dens-cor} to obtain
a quadratic factor $(\B_1,\B_2)$ in $W$ of complexity at most $((32/\delta^4)^{3C_0}, (32/\delta^4)^{2C_0})$
such that the function $g := \E(f|\B_2)$ obeys
$$ |\Lambda_{W,W,W,W}(g,g,g,g) - \Lambda_{W,W,W,W}(\delta, \delta, \delta, \delta)| > \delta^4 / 4.$$
We claim that $\|g\|_{L^\infty(W)} \geq \delta + \frac{\delta}{64}$.  For if this were not the case, then
from Lemma \ref{l1-lambda4} we would have
\[ |\Lambda_{W,W,W,W}(g,g,g,g) - \Lambda_{W,W,W,W}(\delta, \delta, \delta, \delta)| \leq 4 (\delta + \frac{\delta}{64})^3 \|g-\delta\|_{L^1(W)}\]
and hence certainly
$$ \| g - \delta \|_{L^1(W)} \geq \frac{\delta}{32}.$$
Since $g-\delta$ has mean zero, we see (cf. the remarks at the beginning of the proof of Lemma \ref{lininc}) that the maximum value of $g-\delta$ is at least $\frac{\delta}{64}$, a contradiction.  Thus we
indeed have 
\[\| \E(f|\B_2) \|_{L^\infty(W)} \geq \delta + \frac{\delta}{64}.\]
In particular there exists an atom $A$ of $\B_2$ such that $\E_A(f) \geq \E_W(f) + \frac{\delta}{64}$.
The claim now follows from Corollary \ref{linquad}.\endproof

 The above argument was very crude, as it relied on the rather low-technology estimate in Lemma \ref{l1-lambda4}.
In particular, the quadratic structure of $\B_2$ was not used, except in the final step of Corollary \ref{linquad} to
convert the density increment on an atom of $\B_2$ to a density increment on a subspace.  In the next section we refine these computations
by exploiting some ``mixing'' properties of the quadratic factor to obtain some further concentration properties of $g$; only after
obtaining such properties do we invoke the (expensive) Corollary \ref{linquad}.

\section{Quadratic mixing}\label{mixing-sec}

 We have already seen the Chevalley-Warning theorem (Lemma \ref{cwt}) which gives some control on the size of quadratic atoms. It turns out that 
one can do substantially better than this if we assume a non-degeneracy condition on the quadratic phases which define the atom.

 It is convenient to work for now in a homogenized setting, returning to the affine setting later.

\begin{definition}[Homogeneity]  Let $W$ be a vector space (i.e. an affine space with a distinguished origin $0$).  A \emph{homogeneous quadratic phase function} on $W$ is a quadratic phase function $\phi: W \to \R/\Z$ such that $\phi(0) = \nabla \phi(0) = 0$ (thus $\phi(x) = \frac{1}{2} \nabla^2 \phi x \cdot x$).  A \emph{homogeneous linear phase function} on $W$ is a linear phase function $\phi: W \to \R/\Z$ such that $\phi(0) = 0$ (thus $\phi(x) = \nabla \phi \cdot x$).  A \emph{homogenized quadratic factor with complexity $(d_1,d_2)$} on $W$ is any factor which is generated by
$d_1$ homogeneous linear phase functions and $d_2$ homogeneous quadratics, with these $d_1+d_2$ discretized phase functions being linearly independent over $F$.
\end{definition}

 Note that any pure quadratic factor of complexity at most $d$ on a vector space $W$ can be extended to a homogenized quadratic factor of
complexity at most $(d,d)$, simply by taking all the quadratic phases generating the original factor and breaking them up into homogeneous quadratic
and homogeneous linear components (dropping the constant terms, which are not relevant), and then eliminating any linearly dependent terms.

 Now we define the rank of a homogenized quadratic factor.

\begin{definition}[Rank]  Let $\B = \B_{\gamma_1} \vee \ldots \vee \B_{\gamma_{d_1}} \vee \B_{\phi_1} \vee \ldots \vee \B_{\phi_{d_2}}$
be a homogenized quadratic factor of complexity $(d_1,d_2)$ 
on a vector space $W$, generated by $d_1$ homogeneous linear phases $\gamma_i$ and $d_2$ homogeneous quadratic phases $\phi_j$.
We define the \emph{rank} of the factor to be the minimal rank of $\phi$, where $\phi$ ranges over all linear combinations 
\[ \sum_{j = 1}^{d_2} \lambda_j \phi_j,\]
where the $\lambda_j$ are elements of $F$, not all zero,
of $\phi_1,\ldots,\phi_{d_2}$ which are not identically zero. (If $d_2=0$, we define the rank to be infinite.)
\end{definition}

 Intuitively, quadratic factors $\B$ with high rank are highly nonlinear, and as such have a certain amount of ``mixing''. In practise this means that many quantities involving $\B$-measurable functions can be easily understood by working in \emph{configuration space}\footnote{This configuration space can be viewed as a discrete finitary analogue of the $2$-step nilmanifolds which arise naturally in the study of characteristic factors for the $U^3$ norm or four-term recurrence; see \cite{bhk,host-kra,ziegler}.}
 $\T_F^{d_1} \times \T_F^{d_2}$. If $\B = \B_{\gamma_1} \vee \ldots \vee \B_{\gamma_{d_1}} \vee \B_{\phi_1} \vee \ldots \vee \B_{\phi_{d_2}}$ is a homogenized quadratic factor of complexity $(d_1,d_2)$ on a linear space $W$ then we write $\Gamma: W \to \T_F^{d_1}$ and $\Phi: W \to \T_F^{d_2}$ for the maps
$\Gamma(x) := (\gamma_1(x),\ldots,\gamma_{d_1}(x))$ and $\Phi(x) := (\phi_1(x),\ldots,\phi_{d_2}(x))$. If $f : W \rightarrow \C$ is a 1-bounded $\B$-measurable function then we write $\f : \T_F^{d_1} \times \T_F^{d_2} \rightarrow \C$ for the function which satisfies
\[ f(x) = \f(\Gamma(x), \Phi(x))\] for all $x \in W$ and $\f(x_1,x_2) = 0$ if $(x_1,x_2) \neq (\Gamma(x),\Phi(x))$ for any\footnote{The definition of $\f$ outside the range of $\Gamma \times \Phi$ is made merely for definiteness. In later contexts, as the reader may check, $\Gamma \times \Phi : W \to \T_F^{d_1} \times \T_F^{d_2}$ will always be surjective.} $x \in W$. We will adopt this convention of using bold letters to denote functions on configuration space for the rest of the paper without further comment.

 One basic formulation of this principle
is given in Lemma \ref{quad-expect} below. Before we can prove it, we recall a well-known bound on the magnitude of Gauss sums.

\begin{lemma}[Gauss sums]\label{gauss-lem}
Suppose that $W$ is a linear space and that $\phi : W \rightarrow \R/\Z$ is a quadratic form with rank $r$. Then we have the estimate
\[ |\E_x e(\phi(x))| \leq |F|^{-r/2}.\]
\end{lemma} 
\textit{Remark.} Note that the estimate is invariant under adding an arbitrary linear phase to the quadratic form $\phi$.

\textit{Proof.} Write $G_{\phi} := |\E_x e(\phi(x))|$. Squaring and changing variables, we have
\[ G^2_{\phi} = |\E_{x,h} e(\phi(x + h) - \phi(x))|.\]
Using the Taylor expansion \eqref{taylor} and applying the triangle inequality, this gives
\begin{equation}\boxeq G^2_{\phi} \leq \E_x |\E_h e(\nabla \phi(x) \cdot h)| = \P_x(\nabla \phi(x) = 0) = \frac{|\ker(\nabla^2 \phi)|}{|F|^n} = |F|^{-r}.\end{equation}

\begin{lemma}[Expectation on quadratic factors]\label{quad-expect} Let 
$\B$ be a homogenized
quadratic factor of complexity $(d_1,d_2)$ on a linear space $W$, with rank at least $r$.  Let $\Gamma,\Phi$ be the maps from $W$ to configuration space $\T_F^{d_1} \times \T_F^{d_2}$.
Let $f: W \to \C$ be a 1-bounded $\B$-measurable function, and let $\f$ be the corresponding function on configuration space.  Then we have
\[ |\E_W(f) - \E_{\T_F^{d_1} \times \T_F^{d_2}}(\f)| \leq |F|^{(d_1+d_2-r)/2}.\]
\end{lemma}

\begin{proof} We employ a Fourier expansion on configuration space\footnote{Note that $\T_F^d$ is a $d$-dimensional vector space over $F$, and its Pontryagin dual is naturally identified with $F^d$. We persist with the notation $\T_F^{d_1} \times \T_F^{d_2}$ to help the reader, who should remember that any such vector space is being used to label atoms in a quadratic factor. }. Writing
\[ \widehat{\f}(\xi_1,\xi_2) := \E_{(x_1,x_2) \in \T_F^{d_1} \times \T_F^{d_2}} \f (x_1,x_2) e(-\xi_1 \cdot x_1 - \xi_2 \cdot x_2),\]
this allows us to write
\[ f(x) = \sum_{(\xi_1,\xi_2) \in F^{d_1} \times F^{d_2}} 
\widehat \f(\xi_1,\xi_2) e( \xi_1 \cdot \Gamma(x) + \xi_2 \cdot \Phi(x) ).\]
Since $\widehat \f(0,0) = \E_{\T_F^{d_1} \times \T_F^{d_2}}(\f)$, we conclude that
\[ \E_W(f) - \E_{\T_F^{d_1} \times \T_F^{d_2}}(\f) =
\sum_{(\xi_1,\xi_2) \in F^{d_1} \times F^{d_2} \backslash (0,0)} 
\widehat \f(\xi_1,\xi_2) \E_{x \in W} e( \xi_1 \cdot \Gamma(x) + \xi_2 \cdot \Phi(x) ).\]
Now from the rank hypotheses we see that $\xi_1 \cdot \Gamma(x) + \xi_2 \cdot \Phi(x)$ either is a non-constant linear phase,
or is a non-linear quadratic phase of rank at least $r$.  In the former case, the expectation appearing above is zero, whereas
in the latter case the expectation has magnitude at most $|F|^{-r/2}$ by Lemma \ref{gauss-lem}.
Thus by the triangle inequality we have
$$ |\E_W(f) - \E_{\T_F^{d_1} \times \T_F^{d_2}}(\f)| \leq |F|^{-r/2}
\sum_{(\xi_1,\xi_2) \in F^{d_1} \times F^{d_2}} |\widehat \f(\xi_1,\xi_2)|.$$
By Cauchy-Schwarz and Plancherel we have
$$ \sum_{(\xi_1,\xi_2) \in F^{d_1} \times F^{d_2}} |\widehat \f(\xi_1,\xi_2)|
\leq |F|^{(d_1+d_2)/2} \|\f\|_{L^2(\T_F^{d_1} \times \T_F^{d_2})}$$
and the claim now follows from the 1-boundedness of $\f$.
\end{proof}

 We turn now to the somewhat more complicated task of counting 4-term arithmetic progressions using the configuration space, beginning with a heuristic discussion.  
Suppose that $f_0,f_1,f_2,f_3$ are $\B$-measurable functions, and that we wish to compute \[ \Lambda_{W,W,W,W}(f_0,f_1,f_2,f_3).\]  To see what one would expect to get, let us write $f_i(x) = \f_i(\Gamma(x), \Phi(x))$
as before, and expand
$$\Lambda_{W,W,W,W}(f_0,f_1,f_2,f_3) = \E_{x,h \in W} \prod_{i=0}^3 \f_i(\Gamma(x+ih), \Phi(x+ih)).$$
It is then natural to ask what the constraints are between the quantities $\Gamma(x+ih)$ and $\Phi(x+ih)$.  From the linearity of $\Gamma$
and the quadratic nature of $\Phi$ one can easily deduce the constraints
$$ \Gamma(x), \Gamma(x+h), \Gamma(x+2h), \Gamma(x+3h) \hbox{ are in arithmetic progression }$$
and
$$ \Phi(x) - 3\Phi(x+h) + 3 \Phi(x+2h) - \Phi(x+3h) = 0.$$
It turns out that if the rank $r$ is sufficiently large then these are in some sense the ``only'' constraints, and furthermore there is
a certain uniform distribution among all the values of $\Gamma(x+ih)$ and $\Phi(x+ih)$ obeying these constraints.  This leads to
the heuristic formula
\[\Lambda_{W,W,W,W}(f_0,f_1,f_2,f_3) \approx \E_{\begin{array}{l} \scriptstyle x_1, h_1 \in \T_F^{d_1} \\ \scriptstyle
x_{2,0},x_{2,1},x_{2,2},x_{2,3} \in \T_F^{d_2} \\ \scriptstyle x_{2,0}-3x_{2,1}+3x_{2,2}-x_{2,3}=0\end{array}}
\prod_{i=0}^3 \f_i( x_1 + i h_1, x_{2,i} )\]
which can be rearranged using the Fourier transform in the $\T_F^{d_2}$ variables as
\[\begin{split}\Lambda_{W,W,W,W}&(f_0,f_1,f_2,f_3) \approx  \\ &\E_{x_1, h_1 \in \T_F^{d_1}}\sum_{\xi \in F^{d_2}}
\widetilde \f_0(x_1, \xi)
\widetilde \f_1(x_1+h_1, -3\xi)
\widetilde \f_2(x_1+2h_1, 3\xi)
\widetilde \f_3(x_1+3h_1, -\xi)\end{split}\]
where $\widetilde \f$ is the partial Fourier transform of $\f$,
\begin{equation}\label{partial-fourier}
\widetilde \f(x_1,\xi) := \E_{x_2 \in \T_F^{d_2}} \f(x_1,x_2) e(-\xi \cdot x_2).
\end{equation}

 Let us remark that these formulae are closely related to the computations on $2$-step nilmanifolds in \cite{bhk}.  One can view
$\T_F^{d_1} \times \T_F^{d_2}$ as an ``abelian extension'' of the ``Kronecker factor'' $\T_F^{d_1}$, thus creating
a discrete analogue of a $2$-step nilmanifold.  The above formula then is computing $\Lambda$ by taking the Fourier transform in the abelian
extension variable. 

 The next lemma constitutes the rigorous version of the above heuristics.

\begin{lemma}[$\Lambda$ on quadratic factors]\label{quad-expect-4}  Let 
$\B$ be a homogenized
quadratic factor of complexity $(d_1,d_2)$ on a linear space $W$, with rank at least $r$.  Let $\Gamma,\Phi$ be the maps from $W$ to configuration space $\T_F^{d_1} \times \T_F^{d_2}$.  For $i=0,1,2,3$,
let $f_i(x) = \f_i(\Gamma(x), \Phi(x))$ be 1-bounded $\B$-measurable functions.  Then we have
\begin{equation}\label{lambda-approx}
\begin{split}
|\Lambda_{W,W,W,W}&(f_0,f_1,f_2,f_3) - \\
&\E_{x_1, h_1 \in \T_F^{d_1}} \sum_{\xi \in F^{d_2}}\widetilde \f_0(x_1, \xi)
\widetilde \f_1(x_1+h_1, -3\xi)
\widetilde \f_2(x_1+2h_1, 3\xi)
\widetilde \f_3(x_1+3h_1, -\xi)| \\ &\qquad\qquad\qquad\qquad\qquad\qquad\qquad\leq |F|^{(4d_1+4d_2-r)/2},
\end{split}
\end{equation}
where $\widetilde \f_j$ is defined by \eqref{partial-fourier}.
\end{lemma}

\begin{proof}  We use the total Fourier expansion
$$ f_i(x) = \sum_{\lambda_i \in F^{d_1}, \xi_i \in F^{d_2}} \widehat \f_i(\lambda_i, \xi_i) 
e(\lambda_i \cdot \Gamma(x) + \xi_i \cdot \Phi(x))$$
to obtain
\begin{equation}\label{to-use1} \Lambda_{W,W,W,W}(f_0,f_1,f_2,f_3) =
\sum_{\lambda \in (F^{d_1})^4, \xi \in (F^{d_2})^4} m(\lambda,\xi) \prod_{i=0}^3 \widehat \f_i(\lambda_i,\xi_i)\end{equation}
where $\lambda = (\lambda_0,\lambda_1,\lambda_2,\lambda_3)$, $\xi = (\xi_0, \xi_1, \xi_2, \xi_3)$ and
\begin{equation}\label{to-use2}m(\lambda,\xi) := \E_{x,h \in W} e( \sum_{i=0}^3 \lambda_i \cdot \Gamma(x+ih) + \xi_i \cdot \Phi(x+ih) ).\end{equation}
Meanwhile, we may use the Fourier inversion identity
\[ \widetilde{\f}_0(x_1,\xi) = \sum_{\lambda_1} \widehat{\f}_0(\lambda_1,\xi) e(\lambda_1 \cdot x_1) \] together with similar identities for $\widetilde{\f}_1(x_1 + h_1,-3\xi)$, $\widetilde{\f}_2(x_1 + 2h_1,3\xi)$ and $\widetilde{\f}_3(x_1 + 3h_1,-\xi)$
to deduce the formula
\begin{align}\nonumber
 \E_{x_1, h_1 \in \T_F^{d_1}} \sum_{\xi \in F^{d_2}}
\widetilde \f_0(x_1, \xi)
\widetilde \f_1(x_1+h_1, -3\xi)
& \widetilde \f_2(x_1+2h_1, 3\xi)
\widetilde \f_3(x_1+3h_1, -\xi) \\
&=
\sum_{\lambda \in (F^{d_1})^4, \xi \in (F^{d_2})^4} 1_\Sigma(\lambda,\xi) \prod_{i=0}^3 \widehat \f_i(\lambda_i,\xi_i),\label{to-use3}
\end{align}
where $\Sigma \in (F^{d_1})^4 \times (F^{d_2})^4$ is the set of all pairs $(\lambda,\xi)$ such that
\begin{equation}\label{xi-constraints}
3\xi_0 = -\xi_1 = \xi_2 = -3\xi_3 
\end{equation}
and
\begin{equation}\label{lam-constraints}
\lambda_1 + \lambda_2 + \lambda_3 + \lambda_4 = \lambda_2 + 2\lambda_3 + 3\lambda_4 = 0.
\end{equation}
We will shortly show that
\begin{equation}\label{sigma-err}
|m(\lambda,\xi) - 1_{\Sigma}(\lambda,\xi)| \leq |F|^{-r/2}.
\end{equation}
Assuming this, we can compare \eqref{to-use1} with \eqref{to-use3}, bounding the left-hand side of \eqref{lambda-approx} by
$$ |F|^{-r/2} \sum_{\lambda \in (F^{d_1})^4, \xi \in (F^{d_2})^4} \prod_{i=0}^3 |\widehat \f_i(\lambda_i,\xi_i)|.$$
Applying Cauchy-Schwarz and Plancherel as in the proof of the preceding lemma, we can bound
this by $|F|^{(4d_1+4d_2-r)/2}$ as desired.

 It remains to prove \eqref{sigma-err}.  First suppose that \eqref{xi-constraints} fails, so that $1_{\Sigma}(\lambda,\xi) = 0$.  
Then (by a simple inspection) we can find $i' \in \{0,1,2,3\}$ such that
$\sum_{i=0}^3 (i-i')^2 \xi_i \neq 0$.  We can use the change of variables $x = y - i'h$ to write
\[m(\lambda,\xi) := \E_{y,h \in W} e( \sum_{i=0}^3 \lambda_i \cdot \Gamma(y+(i-i')h) + \xi_i \cdot \Phi(y+(i-i')h) ).\]
It then follows from the rank condition that the phase $\sum_{i=0}^3 \xi_i \cdot \Phi(y+(i-i')h)$ contains a non-trivial quadratic
component in $h$ of rank at least $r$. Noting that the linear terms $\lambda_i \cdot \Gamma(y+(i-i')h)$ do not affect the quadratic component
of the phase, we conclude by averaging over $h$ and applying Lemma \ref{gauss-lem} that $m(\lambda,\xi)$ does indeed have magnitude at most $|F|^{-r/2}$.

 Now suppose that \eqref{xi-constraints} holds, but \eqref{lam-constraints} fails, so again $1_{\Sigma}(\lambda,\xi) = 0$.  Then the quadratic nature of $\Phi$
ensures that $\sum_{i=0}^3 \xi_i \cdot \Phi(x+ih) = 0$. Thus 
\begin{eqnarray*} m(\lambda,\xi) & = &\E_{x,h \in W} e(\sum_{i=0}^3 \lambda_i \cdot \Gamma(x+ih)) \\ & = & \E_{x,h \in W}e\big((\lambda_1 + \lambda_2 + \lambda_3 + \lambda_4) \cdot \Gamma(x) + (\lambda_2 + 2\lambda_3 + 3\lambda_4) \cdot \Gamma(h)\big).\end{eqnarray*} The fact that at least one of the vectors $\lambda_1 + \lambda_2 + \lambda_3 + \lambda_4$ and $\lambda_2 + 2\lambda_3 + 3\lambda_4$ is non-zero, combined with the assumed linear independence of the phases $\gamma_1,\dots,\gamma_{d_1}$ comprising $\Gamma$, ensures that $m(\lambda,\xi) = 0$.

 Finally, when \eqref{xi-constraints} and \eqref{lam-constraints} both hold we see that $m(\xi,\lambda) = 1 = 1_{\Sigma}(\lambda,\xi)$ and the claim is trivial.
\end{proof}

 Lemma \ref{quad-expect-4} leads to the following density increment result.

\begin{theorem}[Anomalous AP4-count implies density increment]\label{quad-4-increment} Let
$\B$ be a homogenized quadratic factor of complexity $(d_1,d_2)$ on a linear space $W$ with rank at least $r$.  
Let $f_0,f_1,f_2,f_3: W \to [0,1]$ be 1-bounded non-negative $\B$-measurable functions
which obey the estimates
\begin{equation}\label{fff}
 |\Lambda_{W,W,W,W}(f_0,f_1,f_2,f_3) - \E_W(f_0) \E_W(f_1) \E_W(f_2) \E_W(f_3)| \geq \eta
\end{equation}
and
\begin{equation}\label{gif}
\max_{0 \leq i \leq 3} \E_W(f_i) \leq 6\eta^{1/4}
\end{equation}
for some $\eta$, $1 \geq \eta \geq 2^{40} |F|^{6d_1 + 6d_2 - r}$.  
Then there exists $i$, $0 \leq i \leq 3$, such that one of the following two possibilities hold:

\begin{itemize}

\item \textup{(}medium-sized increment on large subspace\textup{)} There exists an affine subspace $W'$ of $W$ with dimension satisfying
$$ \dim(W') \geq \dim(W) - d_1$$
such that we have the density increment
$$ \E_{W'}(f_i) \geq \E_W(f_i) + 2^{-13}\eta^2.$$

\item \textup{(}large increment on medium-sized subspace\textup{)} There exists a positive integer $K \leq (16/\eta)^3$, and an affine subspace $W'$ of $W$ with 
dimension satisfying
$$ \dim(W') \geq \frac{1}{K+1} \dim(W) - 2(16/\eta)^3 - d_1$$
such that we have the density increment
$$ \E_{W'}(f_i) \geq \E_W(f_i) + 2^{-10} K^{1/3} \eta^{1/4}.$$

\end{itemize}

\end{theorem}

\textit{Remarks.} The constants such as $2^{20}$ appearing here are not best possible.  However, to remove the hypotheses on rank completely will require
an additional argument which we present after proving this theorem.  The density increment obtained here is somewhat better than that
in Proposition \ref{cheap-inc}, for when $K$ is small we do not reduce the dimension of $W$ by as much as in that proposition, and when $K$ 
is large we increase the density on $W$ by significantly more.

\begin{proof}  Let $\Gamma, \Phi$ and $\f_i$ be as before: recall that $\Gamma(x) = (\gamma_1(x),\dots,\gamma_{d_1}(x))$, $\Phi(x) = (\phi_1(x),\dots,\phi_{d_2}(x))$ and that $\f_i$ is a $\B$-measurable function such that \[ \f_i(\Gamma(x),\Phi(x)) = f_i(x).\] Applying Lemma \ref{quad-expect-4} we immediately deduce that
\begin{align}\nonumber
|
\E_{x_1, h_1 \in \T_F^{d_1}} &\sum_{\xi \in F^{d_2}}
\widetilde \f_0(x_1, \xi)
\widetilde \f_1(x_1+h_1, -3\xi)
\widetilde \f_2(x_1+2h_1, 3\xi)
\widetilde \f_3(x_1+3h_1, -\xi)\\
&- \E_W(f_0) \E_W(f_1) \E_W(f_2) \E_W(f_3)| \geq \eta/2.\label{eq8.80}
\end{align}
Now $\widetilde \f_i(x_1,0)$ is close to the average of $f_i$ on the affine space
$\Gamma^{-1}(x_1)$. Indeed applying Lemma \ref{quad-expect} with $f$ replaced by $f1_{\Gamma^{-1}(x_1)}$, so that the corresponding function on configuration space is $\mathbf{f} 1_{x_1 \times \T_F^{d_2}}$, we see that 
\[ |\tilde{\mathbf{f}}_1(x_1,0) - \E_{x \in \Gamma^{-1}(x_1)} f_i(x)| \leq |F|^{(3d_1 + d_2 - r)/2} \leq 2^{-13}\eta^2.\]
Hence if there is some $i$ for which 
\begin{equation}\label{eta-smudge}
 \E_{x_1 \in \T_F^{d_1}} |\widetilde \f_i(x_1,0) - \E_W(f_i)| \geq 2^{-12}\eta^2\end{equation} then
\[ \E_{x_1 \in \T_F^{d_1}} \big| \E_{x \in \Gamma^{-1}(x_1)} f_i(x) - \E_W(f_i)\big| \geq 2^{-13}\eta^2,\]
or in other words
\[ \Vert \E(f_i | \mathcal{B}_{\lin}) - \E_W(f_i) \Vert_{L^1(W)} \geq 2^{-13}\eta^2,\]
where $\mathcal{B}_{\lin}$ is the linear factor of complexity $d_1$ determined by the affine spaces $\Gamma^{-1}(x_1)$, $x_1 \in \T_F^{d_1}$. Lemma \ref{lininc} then tells us that there is some subspace $W'$ with $\dim(W') \geq \dim(W) - d_1$, such that 
\[ \E_{W'}(f_i) \geq \E_W (f_i) + 2^{-13} \eta^2.\]
In this case, then, we have a medium-sized density increment on a large subspace and are done. Suppose, henceforth, that \eqref{eta-smudge} does not hold.

 Now from Lemma \ref{l1-lambda4} we conclude
\begin{align*}
|\E_{x_1, h_1 \in \T_F^{d_1}} 
\widetilde \f_0(x_1, 0)
\widetilde \f_1(x_1+h_1, 0)&
\widetilde \f_2(x_1+2h_1, 0)
\widetilde \f_3(x_1+3h_1, 0) \\
& - \E_W(f_0) \E_W(f_1) \E_W(f_2) \E_W(f_3)| \leq \eta/4
\end{align*}
and hence, from \eqref{eq8.80}, it follows that we must have
\begin{equation}\label{nonzero-mode}
 \E_{x_1, h_1 \in \T_F^{d_1}} \sum_{\xi \in F^{d_2} \backslash 0}
|\widetilde \f_0(x_1, \xi)
\widetilde \f_1(x_1+h_1, -3\xi)
\widetilde \f_2(x_1+2h_1, 3\xi)
\widetilde \f_3(x_1+3h_1, -\xi)| \geq \eta/4.
\end{equation}
From H\"older's inequality and a change of variables we conclude that there exists $i$, 
$0 \leq i \leq 3$, such that
\[ \E_{x_1 \in \T_F^{d_1}} \sum_{\xi \in F^{d_2} \backslash 0} |\widetilde \f_i(x_1,\xi)|^4 \geq \eta/4.\]
It is immediate from Plancherel's identity that
\[ \sum_{\xi \in F^{d_2} \backslash 0} |\widetilde \f_i(x_1,\xi)|^4 \leq 1,\]
and so a simple averaging argument tells us that for a proportion at least $\eta/8$ of the $x_1 \in \T_F^{d_1}$ we have
\begin{equation}\label{eta-8}
\sum_{\xi \in F^{d_2} \backslash 0} |\widetilde \f_i(x_1,\xi)|^4 \geq \eta/8.
\end{equation}
Now we are assuming that \eqref{eta-smudge} does not hold, that is to say
\[ \E_{x_1 \in \T_F^{d_1}} |\widetilde \f_i(x_1,0) - \E_W(f_i)| < 2^{-12}\eta^2.\]
It follows that the proportion of values of $x_1 \in \T_F^{d_1}$ for which it is not the case that 
\begin{equation}\label{eta-flat} 
|\widetilde \f_i(x_1,0) - \E_W(f_i)| \leq 2^{-9}\eta
\end{equation} is less than $\eta/8$.
In particular, there is at least one value of $x_1$ such that both \eqref{eta-8} and \eqref{eta-flat} are satisfied.
Fix this $x_1$, and write $\F_i(x_2) := \f_i(x_1,x_2)$ and $\widehat \F_i(\xi) := \widetilde \f_i(x_1,\xi)$. Note that \eqref{eta-flat} can be written in the form\begin{equation}\label{eta-flat-2}
|\widehat \F_i(0) - \E_W(f_i)| \leq 2^{-12}\eta.
\end{equation}
In particular, in view of \eqref{gif}, we have
\begin{equation}\label{gif-2} |\widehat \F_i(0)| \leq 8 \eta^{1/4}.\end{equation}
At this point we employ some arguments very close to those of Heath-Brown \cite{heath} and Szemer\'edi \cite{szem-3ap}.  From Plancherel's theorem we have
$$ \sum_{\xi \in F^{d_2} \backslash 0} |\widehat \F_i(\xi)|^2 \leq 1.$$
and hence by \eqref{eta-8}
$$ \sum_{\xi \in F^{d_2} \backslash 0: |\widehat \F_i(\xi)|^2 \geq \eta/16} |\widehat \F_i(\xi)|^4 \geq \eta/16.$$
Let us order the $\xi$ in this summation as $\xi_1, \xi_2, \ldots, \xi_J$ in decreasing order of $|\widehat \F_i(\xi)|$, so that
$J \leq 16/\eta$ and
\[ \sum_{j=1}^J |\widehat \F_i(\xi_j)|^4 \geq \eta/16.\]
By the pigeonhole principle and the fact that $\zeta(4/3) \leq 16$, there exists $K$,  $1 \leq K \leq J,$ such that
\[ |\widehat \F_i(\xi_K)| \geq \frac{\eta^{1/4}}{4 K^{1/3} }.\]
Fix this $K$, and set $S := \{ \xi_1,\ldots,\xi_K \}$. We clearly have
\[ \sum_{\xi \in S} |\widehat \F_i(\xi)|^2 \geq K |\widehat \F_i(\xi_K)|^2 \geq  \eta^{1/2} K^{1/3}/16.\]
Thus $S$ has captured a significant amount of $L^2$ energy of $\F_i$ in frequency space; we now look for 
a similar concentration of energy in physical space.
Let $S^\perp \subset \T_F^{d_2}$ be the orthogonal complement of $S$ in $\T_F^{d_2}$, which is a linear subspace of $\T_F^{d_2}$ of codimension
at most $|S|$.  Note that $\{0\} \cup S \subseteq S^{\perp \perp}$.  From the Poisson summation formula and Plancherel we have
$$ \sum_{\xi \in S^{\perp\perp}} |\widehat \F_i(\xi)|^2 = \E_{c \in \T_F^{d_2}} |\E_{x \in c + S^\perp} \F_i(x)|^2$$
and hence
\[ \E_{c \in \T_F^{d_2}} |\E_{x \in c + S^\perp} \F_i(x)|^2 \geq |\widehat \F_i(0)|^2 +  \eta^{1/2} K^{1/3}/16.\]
Now from the positivity of $\F_i$ we have
\[ \E_{c \in \T_F^{d_2}} |\E_{x \in c + S^\perp} \F_i(x)| = \E_{c \in \T_F^{d_2}} \E_{x \in c + S^\perp} \F_i(x) = \widehat \F_i(0),\]
and hence there exists a coset $c+S^\perp$ of $S^\perp$ such that
\[ \E_{x \in c + S^\perp} \F_i(x) \geq \widehat \F_i(0) +  \eta^{1/2} K^{1/3} / 16\widehat \F_i(0).\]
Using \eqref{eta-flat-2} and \eqref{gif-2} it is a simple matter to conclude that
\begin{equation}\label{zoop}
 \E_{x \in c + S^\perp} \F_i(x) \geq \E_W(f_i) +  K^{1/3} \eta^{1/4}/256.
\end{equation}
Now recall that $\F_i(x) = \f_i(x_1,x)$, and that $f_i(x) = \f_i(\Gamma(x),\Phi(x))$. Thus \eqref{zoop} is asserting a density increment for $f_i$ on 
\[ A := \Gamma^{-1}(x_1) \cap \Phi^{-1}(c + S^{\perp}).\]
Note that $A$ is a collection of $\B_S$-atoms, where $\B_S$ is the factor of complexity $(d,K)$ generated by $\gamma_1,\dots,\gamma_d$ and $\xi_1 \cdot \Phi,\dots,\xi_K \cdot \Phi$. To quantify this density increment precisely we must apply Lemma \ref{quad-expect}. Write $g := f_i1_{A}$ and note that if 
\[ \g(t_1,t_2) := 1_{t_1 = x_1} 1_{t_2 \in c + S^{\perp}} \f_i(t_1,t_2)\] then
\[ g(x) = \g(\Gamma(x),\Phi(x)).\]
Applying Lemma \ref{quad-expect} to the function $g$ and noting that
\[ \E_{\T_F^{d_1} \times \T_F^{d_2}}(\g) = |F|^{-d_1} \E \F_i 1_{c + S^{\perp}},\]
we conclude that
\begin{equation}\label{ap-1} ||F|^{-d_1} \E \F_i 1_{c+S^\perp} - \E_W(f_i 1_{A})| \leq |F|^{(d_1+d_2-r)/2}.\end{equation} Applying the same lemma to the function $1_A$, we also have
\begin{equation}\label{ap-2} ||F|^{-d_1} \E 1_{c+S^\perp} -  \E_{ W} 1_A| \leq |F|^{(d_1+d_2-r)/2}.\end{equation}
Now we certainly have 
\begin{equation}\label{est-1}
\E 1_{c+S^\perp} \geq |F|^{-d_2};\end{equation}
this together with \eqref{ap-2} and our assumption on $r$ implies that
\begin{equation}\label{est-2}
\E_W 1_A \geq \textstyle \frac{1}{2} |F|^{-d_2}.
\end{equation}
Combining \eqref{ap-1} with \eqref{est-1} gives
\[ \bigg| \E_{x \in c + S^{\perp}} \F_i(x) - \frac{|F|^{d_1}\E_W (f 1_A)}{\E 1_{c + S^{\perp}}} \bigg | \leq |F|^{(3d_1 + 3d_2 - r)/2},\]
whilst \eqref{ap-2} and \eqref{est-2} together yield
\[ \bigg| \frac{|F|^{d_1}}{\E 1_{c + S^{\perp}}} - \frac{1}{\E_W 1_A} \bigg| \leq \frac{|F|^{(3d_1 + d_2 - r)/2}}{(\E_W 1_A)(\E 1_{c + S^{\perp}})} \leq 2|F|^{(3 d_1 + 5d_2 - r)/2}.\]
Combining these last two inequalities and recalling our assumption on the relation between $\eta$ and $r$, we obtain
\[ |\E_{x \in c+S^\perp} \F_i(x) - \E_{A}( f_i)| \leq 3 |F|^{(3d_1+5d_2-r)/2} \leq 2^{-18}\eta.\]
Inserting this into \eqref{zoop} we conclude that
\[ \E_A(f_i) \geq \E_W(f_i) + 2^{-10} K^{1/3} \eta^{1/4}.\]
In particular there is some atom $A'$ in the factor $\B_S$ such that 
\[ \E_{A'}(f_i) \geq \E_W(f_i) + 2^{-10} K^{1/3} \eta^{1/4}.\]
Applying Corollary \ref{linquad} to the factor $\B_S$, we can then find an affine 
subspace $W'$ of $W$ with dimension satisfying
$$ \dim(W') \geq \frac{1}{K+1} \dim(W) - \frac{2K+d_1}{K+1} - K \geq \frac{1}{K+1} \dim(W) - 2(16/\eta)^3 - d_1$$
such that
$$ \E_{W'}(f_i) \geq \E_W(f_i) + 2^{-10} K^{1/3} \eta^{1/4}.$$
The claim follows.
\end{proof}

 The last result was proved under two assumptions, namely homogeneity and a rank condition. We now remove these hypotheseses. To remove the latter we first need a lemma.

\begin{lemma}[Rank lemma]\label{rank-lemma}  Let $G$ be a linear space, and let $M_1,\ldots,M_d: G \to \widehat G$ be self-adjoint linear
transformations.  Let $r \geq 0$ be arbitrary.  Then there exists a linear subspace $H$ of $G$ of codimension at most $rd$, and
self-adjoint linear transformations $M'_1, \ldots, M'_{d'}: H \to \widehat H$ for some $0 \leq d' \leq d$, such that we have
the rank condition
\begin{equation}\label{rank-cond}
\rank( a_1 M'_1 + \ldots + a_{d'} M'_{d'} ) > r
\end{equation}
whenever $(a_1,\ldots,a_{d'}) \in F^{d'} \setminus \{0\}$. In particular $M'_1,\dots,M'_d$ are linearly independent. Moreover, interpreting $M_i|_H$ as a map into $\hat{H}$, we have $M_i|_H \in \langle M'_1,\dots,M'_d\rangle$ for each $i$.
\end{lemma}

\begin{proof} We induct on $d$.  The case $d=0$ is vacuously true, so suppose $d \geq 0$ and the claim has already been proven for $d-1$.
We may assume that
$$  \rank( a_1 M_1 + \ldots + a_{d} M_{d} ) \leq r$$
for some $a_1,\ldots,a_d$ not all zero, since otherwise we could just set $d'=d$, $H=G$, and $M'_j=M_j$.  By symmetry and scaling we can take
$a_d=1$.  If we then let $G'$ be the kernel of $a_1 M_1 + \ldots + a_d M_d$ then $G'$ has codimension at most $r$, and when restricted to $G'$
the transformation $M_d$ is a linear combination of the $M_1,\ldots,M_{d-1}$ and can thus be safely omitted.  The claim then follows from applying the induction hypothesis to $G'$ and $M_1,\ldots,M_{d-1}$.
\end{proof}

\begin{theorem}[Anomalous AP4-count implies density increment, II]\label{quad-4-increment-2} Let $W_0,\! W_1$, $W_2$, $W_3$
be a progression of affine spaces, and on each $W_i$ let $\B_i$ be a pure quadratic factor of complexity at most $d$,
and let $f_i: W_i \to \R^+$ be 1-bounded non-negative $\B_i$-measurable functions
which obey the estimates
\begin{equation}\label{fff-alt}
 |\Lambda_{W_0,W_1,W_2,W_3}(f_0,f_1,f_2,f_3) - \E_{W_0}(f_0) \E_{W_1}(f_1) \E_{W_2}(f_2) \E_{W_3}(f_3)| \geq \eta
\end{equation}
and
\begin{equation}\label{gif-alt}
\max_{0 \leq i \leq 3} \E_{W_i}(f_i) \leq 3\eta^{1/4}
\end{equation}
for some $\eta$, $0 < \eta \leq 1$.
Then there exists $i$, $0 \leq i \leq 3$, such that one of the following two alternatives holds:

\begin{itemize}

\item \textup{(}medium-sized increment on large subspace\textup{)} There exists an affine subspace $W'_i$ of $W_i$ with dimension satisfying
\[ \dim(W'_i) \geq \dim(W_i) - 200 d^2 - 8d \log_{|F|} \frac{2^{42}}{\eta}\]
such that we have the density increment
\[ \E_{W'_i}(f_i) \geq \E_{W_i}(f_i) + 2^{-23}\eta^3.\]

\item \textup{(}Large increment on a medium-sized subspace\textup{)} There exists a positive integer $K \leq (64/\eta)^3$, and an affine subspace $W'_i$ of $W_i$ with dimension satisfying 
$$ \dim(W'_i) \geq \frac{1}{K+1} \dim(W_i) - 2(64/\eta)^3 - 200d^2 - 8 d \log_{|F|} \frac{2^{42}}{\eta}$$
such that we have the density increment
$$ \E_{W_i'}(f_i) \geq \E_{W_i}(f_i) + 2^{-13} K^{1/3} \eta^{1/4}.$$

\end{itemize}

\end{theorem}

\begin{proof}  By translating each of the $W_i$ (and also $f_i$ and $\B_i$) we may assume that $W_0=W_1=W_2=W_3=G$ is a vector space.
Let $\B = \B_0 \vee \B_1 \vee \B_2 \vee \B_3$, so that $\B$ is a quadratic factor generated by $4d$ quadratic phases
$\phi_1,\ldots,\phi_{4d}$ (adding dummy phases if necessary).
Let $r$ be the integer part of $50d +  \log_{|F|} \frac{2^{42}}{\eta}$.
We use Lemma \ref{rank-lemma} to find a linear subspace $H$ of $G$ of dimension
at least 
$$\dim(H) \geq \dim(G) - 4dr \geq \dim(G) - 200 d^2 - 4d \log_{|F|} \frac{2^{42}}{\eta}$$ 
and
self-adjoint matrices $M'_1,\ldots,M'_{d'}$ for some $d'$,  $0 \leq d' \leq 4d$, obeying the rank condition \eqref{rank-cond}, such
that each of the Hessians $\nabla^2 \phi_1,\ldots, \nabla^2 \phi_{4d}$ when restricted to $H$ becomes a linear combination of the
$M'_1,\ldots,M'_{d'}$.

 Let $\B'$ be the linear factor generated by the cosets of $H$. We may assume that
\begin{equation}\label{effigy}
\| \E(f_i|\B') - \E_G(f_i) \|_{L^1(G)} \leq 2^{-22}\eta^3;
\end{equation}
for each $i$: if not then Lemma \ref{lininc} provides the claimed medium-sized density increment on a large subspace and we are done.

 Assuming then that \eqref{effigy} holds, it follows from Lemma \ref{l1-lambda4} that
$$ |\Lambda_{G,G,G,G}(\E(f_0|\B'),\E(f_1|\B'),\E(f_2|\B'),\E(f_3|\B')) - \E_G(f_0) \E_G(f_1) \E_G(f_2) \E_G(f_3)| \leq \eta/2,$$
and hence by \eqref{fff-alt} that
$$ |\Lambda_{G,G,G,G}(f_0,f_1,f_2,f_3) - \Lambda_{G,G,G,G}(\E(f_0|\B'),\E(f_1|\B'),\E(f_2|\B'),\E(f_3|\B'))| \geq \eta/2.$$
We can rewrite the left-hand side here as
\begin{align*} \big|\E_{x,h \in G/H} \big(\Lambda_{x+H,x+h+H,x+2h+H,x+3h+H}&(f_0,f_1,f_2,f_3) \\ &- \E_{x+H}(f_0) \E_{x+h+H}(f_1) \E_{x+2h+H}(f_2) \E_{x+3h+H}(f_3)\big) \big|,\end{align*}
so for a proportion at least $\eta/4$ of the pairs $(x,h)$ in $G/H \times G/H$ we have
\begin{align}\nonumber
 |\Lambda_{x+H,x+h+H,x+2h+H,x+3h+H}&(f_0,f_1,f_2,f_3) \\ & - \E_{x+H}(f_0) \E_{x+h+H}(f_1) \E_{x+2h+H}(f_2) \E_{x+3h+H}(f_3)| \geq \eta/4.\label{progression}
\end{align}
Now there must be a pair $(x,h)$ satisfying \eqref{progression}, and such that
\begin{equation}\label{827star} |\E_{x+ih+H}(f_i) - \E_G(f_i)| \leq 2^{-18}\eta^2 \hbox{ for } i=0,1,2,3.\end{equation}
Indeed if there were not then for some $i$ we would have
\[ |\E_{x+ih+H}(f_i) - \E_G(f_i)| > 2^{-18}\eta^2\] for a proportion at least $\eta/16$ of the pairs $(x,h) \in G/H \times G/H$. This would be contrary to \eqref{effigy}.
We now take advantage of the affine invariance of our setup by translating each of the $f_i$ so that $x=h=0$, so that \eqref{progression} becomes
$$  |\Lambda_{H,H,H,H}(f_0,f_1,f_2,f_3) - \E_{H}(f_0) \E_{H}(f_1) \E_{H}(f_2) \E_{H}(f_3)| \geq \eta/4.$$
Now observe from \eqref{taylor} that each of the phase functions $\phi_1,\ldots,\phi_{4d}$, when restricted to $H$, is equal to a linear 
combination of the homogeneous quadratic phases $M'_1 x \cdot x, \ldots, M'_{d'} x \cdot x$ plus a homogeneous linear phase, plus a constant.  By
collecting all these homogeneous quadratic and linear phases together, and omitting any which are linearly dependent (note that the rank condition on the $M'_i$ ensure that the quadratic phases have no such linear dependence) we can thus find a homogenized quadratic 
factor $\dot \B$ on $H$ of complexity at most $(4d,d')$ which has rank greater than $r$, such that all the
phases $\phi_1,\ldots,\phi_{4d}$ are $\dot\B$-measurable on $H$. In particular $f_0,f_1,f_2,f_3$ are also $\dot \B$-measurable.
The claim now follows from Theorem \ref{quad-4-increment} (with $W$ replaced by $H$, and $\eta$ replaced by $\eta/4$) and \eqref{827star}.
\end{proof}

 We are finally able to prove Theorem \ref{dens-inc} (and thus, by the argument at the start of \S \ref{lambda-sec}, Theorem \ref{main}).

\textit{Proof of Theorem \ref{dens-inc}.}
We apply Corollary \ref{dens-cor} to obtain
a quadratic factor $(\B_1,\B_2)$ in $W$ of complexity at most $(d_1,d_2) := ((32/\delta^4)^{3C_0}, (32/\delta^4)^{2C_0})$
such that the function $g := \E(f|\B_2)$ obeys
$$ |\Lambda_{W,W,W,W}(g,g,g,g) - \Lambda_{W,W,W,W}(\delta, \delta, \delta, \delta)| > \delta^4 / 4.$$
We may assume that
\begin{equation}\label{b1-razor}
\| \E(f|\B_1) - \delta\|_{L^1(W)} \leq 2^{-42}\delta^{16},
\end{equation}
since otherwise the claim would follow from Lemma \ref{lininc}.  In particular we see from Lemma \ref{l1-lambda4} that
$$ |\Lambda_{W,W,W,W}(\E(f|\B_1),\E(f|\B_1),\E(f|\B_1),\E(f|\B_1)) - \Lambda_{W,W,W,W}(\delta, \delta, \delta, \delta)| \leq \delta^4 / 8$$
and thus
$$ |\Lambda_{W,W,W,W}(g,g,g,g) - \Lambda_{W,W,W,W}(\E(f|\B_1), \E(f|\B_1), \E(f|\B_1), \E(f|\B_1))| > \delta^4 / 8.$$
We can rewrite the left-hand side as
$$ |\E_{W_0,W_1,W_2,W_3} (\Lambda_{W_0,W_1,W_2,W_3}(g,g,g,g) - \E_{W_0}(f) \E_{W_1}(f) \E_{W_2}(f) \E_{W_3}(f))|,$$
where the $W_0,W_1,W_2,W_3$ range over all quadruples of atoms of $\B_1$ in arithmetic progression.  For a proportion at least $\delta^4/16$ of these
progressions we have
\[ |\Lambda_{W_0,W_1,W_2,W_3}(g,g,g,g) - \E_{W_0}(f) \E_{W_1}(f) \E_{W_2}(f) \E_{W_3}(f)| \geq \delta^4/16.\] By a simple averaging argument using \eqref{b1-razor} we can find a progression $W_0,W_1,W_2,W_3$ with this property such that
$$ |\E_{W_i}(f) - \E_W(f)| \leq 2^{-36}\delta^{12}$$
for all $i=0,1,2,3$.  The claim now follows from Theorem \ref{quad-4-increment-2} with $\eta := \delta^4/16$.\endproof

\section{Acknowledgement} The authors would like to thank Roger Heath-Brown for a careful reading of the paper.

\end{document}